\documentclass[a4paper]{article}
\usepackage{amsmath,amsthm,amssymb}

\usepackage{color}
\usepackage[dvips]{graphicx} 
\usepackage{version}
\newenvironment{NB}{
\color{red}{\bf NB}. \footnotesize
}{}
\excludeversion{NB}

\usepackage[all]{xy}

\newtheorem{thm}{Theorem}[section]
\newtheorem{defn}[thm]{Definition}
\newtheorem{ex}[thm]{Example}
\newtheorem{prop}[thm]{Proposition}
\newtheorem{cor}[thm]{Corollary}
\newtheorem{lem}[thm]{Lemma}

\newtheorem{rem}[thm]{Remark}

\newcommand{\mr}[1]{{\mathrm{#1}}}

\newcommand{\bb}[1]{{\mathbb{#1}}}

\newcommand{\Z}{\bb{Z}}
\newcommand{\Zh}{\bb{Z}_{\mathrm{h}}}
\newcommand{\C}{\bb{C}}

\newcommand{\R}{\bb{R}}

\newcommand{\vr}{v_{\mathrm{r}}}
\newcommand{\vl}{v_{\mathrm{l}}}

\newcommand{\Ih}{I_{\mathrm{h}}}
\newcommand{\eh}{e_{\mathrm{h}}}
\newcommand{\es}{e_{\mathrm{s}}}

\newcommand{\tE}{\tilde{\mathcal{E}}}
\newcommand{\plt}{\sigma,\lambda,\theta}
\newcommand{\plmt}{\sigma,\lambda,\nu,\theta}

\newcommand{\tsig}{\tilde{\sigma}}
\newcommand{\tlam}{\tilde{\lambda}}

\usepackage{pst-all}

\numberwithin{equation}{section}

\newcommand{\vlam}{\vec{\lambda}}

\title{Refined open non-commutative Donaldson-Thomas invariants for small crepant resolutions}
\author{Kentaro Nagao\\
RIMS, Kyoto University\\
  Kyoto 606-8502, Japan
}

\begin{document}

\maketitle

\begin{abstract}
The aim of this paper is to study analogs of noncommutative Donaldson-Thomas invariants corresponding to the refined topological vertex for small crepant resolutions of toric Calabi-Yau $3$-folds.
We give three definitions of the invariants which are equivalent to each others and provide ``wall-crossing" formulas for the invariants.
In particular, we get normalized generating functions which are unchanged under ``wall-crossing".
\end{abstract}

\section*{Introduction}
Donaldson-Thomas theory (\cite{thomas-dt}) is intersection theory on the moduli spaces of ideal sheaves on a smooth variety, which is conjecturally equivalent to Gromov-Witten theory (\cite{mnop}). 
For a Calabi-Yau $3$-fold, the virtual dimension of the moduli space is zero and hence Donaldson-Thomas invariants are said to be 
counting invariants of ideal sheaves. 
It is known that they coincide with the weighted Euler characteristics of the moduli spaces weighted by the Behrend's functions (\cite{behrend-dt}).  
Recently, study of Donaldson-Thomas invariants of Calabi-Yau $3$-folds using categorical methods has been being developed (\cite{joyce-4,joyce-3cy,toda-limit,toda-rational,ks,joyce-song} etc.).

On the other hand, a smooth variety $Y$ sometimes has a noncommutative associative algebra $A$ such that the derived category of coherent sheaves on $Y$ is equivalent to the derived category of $A$-modules. 
Derived McKay correspondence (\cite{kapranov-vasserot,bkr}) and Van den Bergh's noncommutative crepant resolutions (\cite{vandenbergh-nccr}) are typical examples. 
In such cases, B. Szendroi proposed to study counting invariants of $A$-modules ({\it noncommutative Donaldson-Thomas invariants}) and relations with the original Donaldson-Thomas invariants on $Y$ (\cite{szendroi-ncdt}). 
In \cite{nagao-nakajima,3tcy}, we provided wall-crossing formulas which relate generating functions of Donaldson-Thomas and noncommutative Donaldson-Thomas invariants for small crepant resolutions of toric Calabi-Yau $3$-folds. (We say a resolution of a $3$-fold is {\it small} if the dimension of each fiber is less than or equal to $1$.)

The aim of this paper is to propose new invariants generalizing noncommutative Donaldson-Thomas invariants and to provide ``wall-crossing formulas" for small crepant resolutions of toric Calabi-Yau $3$-folds. 
We have two directions of generalizations:
\begin{itemize}
\item ``open" version\footnote{The word ``open'' stems from such terminologies as ``open topological string theory''. According to \cite{TV}, open topological string partition function is given by summing up the generating functions of these invariants over Young diagrams.} : corresponding to counting invariants of sheaves on $Y$  with non-compact supports\footnote{As far as the author knows, there is no definition of ``open" invariants for general Calabi-Yau $3$-folds.},
\item refined version : corresponding to refined topological vertex (\cite{RTV})\footnote{See \cite{behrend_bryan_szendroi} for a geometric definition of refined invariants.}.
\end{itemize}

Let $Y\to X$ be a projective small crepant resolution of an affine toric Calabi-Yau $3$-fold. 
Recall that giving an affine toric Calabi-Yau $3$-fold is equivalent to giving a convex lattice polygon. 
Existence of a small crepant resolution is equivalent to absence of interior lattice points in the polygon. 
It is easy to classify such polygons and $X$ is one of the following:
\begin{itemize}
\item $X=X_{L^+,L^-}:=\{\mr{x}\mr{y}=\mr{z}^{L^+}\mr{w}^{L^-}\}\subset \C^4$ for $L^+>0$ and $L^-\geq 0$, or 
\item $X=X_{(\Z/2\Z)^2}:=\C^3/(\Z/2\Z)^2$ where $(\Z/2\Z)^2$ acts on $\C^3$ with weights $(1,0)$, $(0,1)$ and $(1,1)$.
\end{itemize}
\begin{figure}[htbp]
  \centering
  \input{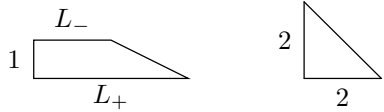}
  \caption{Polygons for $X_{L^+,L^-}$ and $X_{(\Z/2\Z)^2}$}
\label{}
\end{figure}
In this paper, we study the first case.
We put $L:=L^++L^-$.
Note that $X_{1,1}$ is called the {\it conifold} and $X_{L,0}$ is isomorphic to $\C\times\C^2/(\Z/L\Z)$.

Given a pair of Young diagrams ${\nu}=(\nu_+,\nu_-)$ and an $L$-tuple of Young diagrams 
\[
{\lambda}=\Bigl(\lambda^{(1/2)},\ldots,\lambda^{(L-1/2)}\Bigr),
\]
the generating function of {\it refined open noncommutative Donaldson-Thomas invariants} (roncDT, in short)
\[
\mathcal{Z}^Y_{\lambda,{\nu}}\bigl(\vec{q}\,\bigr)
=\mathcal{Z}^Y_{\lambda,{\nu}}\bigl(q_+,q_-,q_1\ldots,q_{L-1}),
\]
which is denoted by $\mathcal{Z}^{\mathrm{RTV}}_{\sigma,\lambda,\nu}$ in the body of this paper, is defined by counting the number of the following data:
\begin{itemize}
\item an $(L-1)$-tuple of Young diagrams $\vec{\nu}=(\nu^{(1)},\ldots,\nu^{(L-1)})$, and 
\item an $L$-tuple of $3$-dimensional Young diagrams $\vec{\Lambda}=\bigl(\Lambda^{(1/2)},\ldots,\Lambda^{(L-1/2)}\bigr)$ such that $\Lambda^{(j)}$ is of type $\bigl(\lambda^{(j)},\nu^{(j+1/2)},{}^{\mathrm{t}}\nu^{(j-1/2)}\bigr)$ or $\bigl(\lambda^{(j)},{}^{\mathrm{t}}\nu^{(j-1/2)},\nu^{(j+1/2)}\bigl)$ (see \S \ref{subsub_515} for the details).
\end{itemize}
Such data parametrize torus fixed ideal sheaves on the small crepant resolution $Y$. 
In particular, 
\[
\mathcal{Z}^Y_{{\emptyset},{\emptyset}}\bigl(\vec{q}\,\bigr)\big|_{q_+=q_-}
\]
coincides with the generating function of Euler characteristic versions of the Donaldson-Thomas invariants of $Y$\footnote{The Euler characteristic version of Donaldson-Thomas invariant coincides with the Donaldson-Thomas invariant up to sign (\cite{mnop}).}.

Let $A$ be a noncommutative crepant resolution of $X$ and $\Zh$ denote the set of half integers and $\theta\colon \Zh\to\Zh$ be a bijection such that $\theta(h+L)=\theta(h)+L$ and such that
\[
\theta(1/2)+\cdots+\theta(L-1/2)=1/2+\cdots+(L-1/2).
\]
In this paper, we will define generating functions $\mathcal{Z}^A_{{\lambda},{\nu},\theta}\bigl(\vec{q}\,\bigr)$, which are denoted by $\mathcal{Z}_{\sigma,\lambda,\nu,\theta}\bigl(\vec{q}\bigr)$ in the body of this paper (see \S \ref{subsec_gen_func}), 
such that 
\begin{itemize}
\item $\mathcal{Z}^A_{{\emptyset},{\emptyset},\mathrm{id}}\bigl(\vec{q}\,\bigr)\big|_{q_+=q_-=q_0^{1/2}}$ coincides with the generating function $\mathcal{Z}^{A}_{\mathrm{NCDT},\mathrm{eu}}$ of Euler characteristic versions\footnote{The Euler characteristic version of noncommutative Donaldson-Thomas invariant coincides with the noncommutative Donaldson-Thomas invariant up to sign (\cite{3tcy, ncdt-brane}).} of noncommutative Donaldson-Thomas invariants for the noncommutative crepant resolution $A$ (see \cite{ncdt-brane} and Remark \ref{rem_brane}), 
\item $``\displaystyle{\lim_{\theta\to\infty}}"\,\mathcal{Z}^A_{{\lambda},{\nu},\theta}\bigl(\vec{q}\,\bigr)=\mathcal{Z}^Y_{{\lambda},{\nu}}\bigl(\vec{q}\,\bigr)$ (Theorem \ref{thm_limit})\footnote{The limit in this equation is, in fact, equivalent to a limit in the space of stability conditions for the category of finite dimensional $A$-modules.}\footnote{A moduli space of stalbe $A$-modules with the specific numerical data gives a crepant resolution of $X$(\cite{ishii-ueda}). The direnction in which we take limit in the space of stability conditions determines a stability paremeter in the construction of a crepant resolution.}.
\end{itemize}
Moreover, for $i\in I:=\Z/L\Z$ we can define the new bijection $\mu_i(\theta)\colon \Zh\to\Zh$ (see \S \ref{subsub_121}) and 
\begin{itemize}
\item $\mathcal{Z}^A_{{\lambda},{\nu},\mu_i(\theta)}\bigl(\vec{q}\,\bigr)\big/\mathcal{Z}^A_{{\lambda},{\nu},\theta}\bigl(\vec{q}\,\bigr)$ is given explicitly (Theorem \ref{thm_hex} and \ref{thm_quad}). 
\end{itemize}
In \cite{nagao-nakajima,3tcy}, we realized $\mathcal{Z}^A_{{\emptyset},{\emptyset},\theta}\bigl(\vec{q}\,\bigr)\big|_{q_+=q_-}$ as generating functions of virtual counting of certain moduli spaces and these moduli spaces are constructed using geometric invariant theory.
In this story, $\theta$ determines a chamber in the space of stability parameters and the chamber corresponding to $\theta$ is adjacent to the chamber corresponding to $\mu_i(\theta)$ by a single wall.
This is the reason we call Theorem \ref{thm_hex} and \ref{thm_quad} as ``wall-crossing" formulas even though our definition of the invariants and the proof of the formula are given in combinatorial ways. 
In fact, in the subsequent paper \cite{NCDTviaVO} we provide an alternative geometric definition, in which $\theta$ determines a chamber in the space of Bridgeland's stability conditions for the category of finite dimensional $A$-modules.

As a corollary of the ``wall-crossing" formula, we get
\begin{itemize}
\item 
$\mathcal{Z}^A_{{\lambda},{\nu},\theta}
\,\Big/
\mathcal{Z}^A_{{\lambda},{\emptyset},\theta}
=
\mathcal{Z}^Y_{{\lambda},{\nu}}
\,\Big/
\mathcal{Z}^Y_{{\lambda},{\emptyset}}
$ \quad (Corollary \ref{cor_413} and \ref{cor_56})
\end{itemize}
for any $\theta$, ${\lambda}$ and ${\nu}$, and
\begin{itemize}
\item
$\Bigl(\mathcal{Z}^A_{{\lambda},{\nu},\theta}
\,\Big/
\mathcal{Z}^A_{{\emptyset},{\emptyset},\theta}\Bigr)\Big|_{q_+=q_-}
=
\Bigl(\mathcal{Z}^Y_{{\lambda},{\nu}}
\,\Big/
\mathcal{Z}^Y_{{\emptyset},{\emptyset}}\Bigr)\Big|_{q_+=q_-}
$ \quad (Corollary \ref{cor_415} and \ref{cor_58})
\end{itemize}
for any $\theta$, $\lambda$ and ${\nu}$ such that $c_\lambda[j]=0$ for any $j$ (see \S \ref{subsub_core} for the notation).
According to the result in \cite{nagao-nakajima,3tcy}, these formulas should be interpreted as stability of the normalized generating functions under ``wall-crossing".
We can find such stability of normalized generating functions in other contexts such as flop invariance and DT-PT correspondence. 
Categorical interpretations of such normalized generating functions and their stability are waited.

\smallskip
Now, we summarize the prior study on noncommutative Donaldson-Thomas invariants : 
\begin{itemize}
\item
Szendroi's formula on the generating function of noncommutative Donaldson-Thomas invariants of the conifold was shown by B. Young in a purely combinatorial way (\cite{young-conifold}). The main tool is an operation called {\it dimer shuffling}.
\item
In \cite{young-mckay}, the authors generalized Szendroi-Young formula for $X_{L,0}$ and $X_{(\Z/2\Z)^2}$. 
The method is different from the one used in \cite{young-conifold} : they use vertex operator method.
\item
In \cite{nagao-nakajima}, we gave an interpretation of Szendroi-Young formula as a consequence of the wall-crossing formula.
From our point of view, the argument in \cite{young-conifold} is obtained by translating the argument in \cite{nagao-nakajima} into combinatorial language by localization.
In particular, dimer shuffling is nothing but ``mutation" in the categorical language.
\item
In \cite{3tcy}, we generalized the results in \cite{nagao-nakajima} for arbitrary small crepant resolutions of toric Calabi-Yau $3$-folds.
\item
In \cite{joyce-song}, the authors study noncommutative Donaldson-Thomas invariants of small crepant resolutions of toric Calabi-Yau $3$-folds as examples of their theory of generalized Donaldson-Thomas invariants.
\item
In \cite{dimofte-gukov}, T. Dimofte and S. Gukov provided a refined version of Szendroi-Young formula for the conifold.
\item See \cite{jafferis-moore,chuang-jafferis,AOVY,CP,aganagic_yamazaki,dimofte_gukov_soibelman} for developments in physics. 
\end{itemize}
In this paper, we define the roncDT invariants using a dimer model (\S \ref{sec_dimer_model}), which is purely combinatorial.

\begin{NB}
In this paper, we give the following three definitions of the invariants:
\begin{itemize}
\item 
In \S \ref{sec_dimer_model}, we give the definition using a dimer model. This is purely combinatorial.
\item 
For $\sigma,\lambda,\nu,\theta$, a module $M_{\sigma,\lambda,\nu,\theta}$ of a noncommutative crepant resolution is given. 
In \S \ref{sec_moduli}, we construct moduli spaces of finite dimensional quotient modules of $M_{\sigma,\lambda,\nu,\theta}$. 
Taking the Euler characteristic of these moduli spaces is the second definition.
\item 
For $\sigma,\lambda,\nu$, we define a ``framed " quiver $Q_{\sigma,\lambda,\nu}$ given by adding one extra vertex and some extra edges from that vertex to the noncommutative crepant resolution.
Regarding $\theta$ as a parameter of stability, we construct moduli spaces of finite dimensional modules of this quiver with a one dimensional vector space on the extra vertex.
Taking the Euler characteristic of these moduli spaces is the third definition.
\end{itemize}
\end{NB}
In \S \ref{sec_crystal}, we give an interpretation of the dimer model as a crystal melting model.
We construct an $A$-module $M^{\mr{max}}_{\sigma,\lambda,\nu,\theta}$ such that giving a dimer configuration is equivalent to giving a finite dimensional torus invariant quotient module of $M^{\mr{max}}_{\sigma,\lambda,\nu,\theta}$.
Hence the roncDT invariant coincides with the Euler characteristic of the moduli space of finite dimensional quotient modules of $M^{\mr{max}}_{\sigma,\lambda,\nu,\theta}$ (see \cite{NCDTviaVO})\footnote{From the geometric point of view, the crystal melting model is more natural.
But in this paper we adapt the definition using the dimer model since it is more convenient when we prove some technical lemmas, which we also use in \cite{NCDTviaVO}.}\footnote{In the case when $\nu_+=\nu_-=\emptyset$, the moduli spaces have symmetric obstruction theory and the invariant in this paper coincides with the weighted Euler characteristic up to sign.}.

\begin{NB}
In \S \ref{sec_moduli}, we construct moduli spaces of finite dimensional quotient modules. 
The torus fixed point set of the moduli spaces are parametrized by dimer configurations and hence the generating function defined by the dimer model coincides with that of the Euler characteristic of the moduli spaces. 
Moreover, in particular cases, we construct symmetric obstruction theory and the generating function by the dimer model coincides with that of virtual counting of the moduli spaces, which are more sophisticated invariants.
\end{NB}%

In \S \ref{sec_shuffling_wc}, we introduce the notion of dimer shuffling to prove the ``wall-crossing" formula (Theorem \ref{thm_hex} and \ref{thm_quad}),
which is the first main result of this paper.

Finally we study the limit behavior of the dimer model in \S \ref{sec_rtv}.
The second main result is that the generating function given by the refined topological vertex for $Y$ appears as the limit (Theorem \ref{thm_limit}).


While preparing the papers, the author was informed from J. Bryan that he and his collaborators C. Cadman and B. Young provided an explicit formula of 
$\mathcal{Z}^A_{{\lambda},{\nu},\mathrm{id}}\big|_{q_+=q_-}$ for $X_{L,0}$ and $X_{(\Z/2\Z)^2}$ using vertex operator methods (\cite{bryan-young}).
In the subsequent paper \cite{NCDTviaVO}, we provide an explicit formula of $\mathcal{Z}^A_{{\lambda},{\nu},\theta}$ for $X_{L_+,L_-}$ using vertex operator methods.

A physicist may refer to \cite{nagao-yamazaki}, in which we explain the result of this paper in a physical context.


\begin{NB}
The is a subsequent paper to \cite{open_conifold}. 

Let $X$ be a $3$-fold with Gorenstein singularities with the trivial canonical bundle.
A crepant resolution $Y\to X$ is said to be {\it small} if the dimensions of the fibers are less or equal than $1$. 
Assume $X$ has a small crepant resolution, then $X$ has a noncommutative crepant resolution $A$ (\cite{vandenbergh-nccr}). 
In particular, the derived category of coherent sheaves on $Y$ is equivalent to the derived category of $A$-modules.

Donaldson-Thomas theory for a Calabi-Yau $3$-fold $Y$ is a counting theory of coherent sheaves on $Y$.
Recently, study of Donaldson-Thomas invariants of Calabi-Yau $3$-folds has been developing from categorical point of view by D. Joyce, Y. Toda, Kontsevich-Soibelman,,, (\cite{joyce-4,joyce-3cy,toda-limit,toda-rational,ks,joyce-song-1,joyce-song-2,joyce-song-3} etc.). 
In the cases as above, B. Szendroi proposed to study counting invariants of $A$-modules ({\it noncommutative Donaldson-Thomas invariants}) and relations with the original Donaldson-Thomas invariants on $Y$ (\cite{szendroi-ncdt}). 
In \cite{nagao-nakajima,3tcy} we provided a wall-crossing formula which relates generating functions of Donaldson-Thomas and noncommutative Donaldson-Thomas invariants.

Let $Y\to X$ be a projective small crepant resolution of toric Calabi-Yau $3$-fold. 
Then $X$ is one of the following:
\begin{itemize}
\item $X=X_{a,b}:=\{xy=z^aw^b\}\subset \C^4$ for $a\geq 0$ and $b>0$, or 
\item $X=X_{\Z/2\Z}:=\C^3/(\Z/2\Z)$ where $\Z/2\Z$ acts on $\C^3$ with weights $(1,0)$, $(0,1)$ and $(1,1)$.
\end{itemize}
In the case $(a,b)=(1,1)$, $X_{1,1}$ is called the {\it conifold}.
In \cite{open_conifold}, 
we defined generating functions
\[
\mathcal{Z}_{n,{\lambda}}(q_+,q_-,q_1)
\]
for a nonnegative integer $n$ and  a quadruplet $\vlam$ of Young diagrams so that 
\begin{itemize}
\item $\mathcal{Z}_{0,{\emptyset}}(q_+.q_-,q_1)|_{q_+=q_-=q_0^{1/2}}$ coincides with the generating function of noncommutative Donaldson-Thomas invariants for the conifold, 
\item $``\displaystyle{\lim_{n\to\infty}}"\,\mathcal{Z}_{n,{\lambda}}(q_+,q_-,q_1)=\mathcal{Z}_{{\lambda}}(q_+,q_-,q_1)$, where the right-hand side is given by the refined topological vertex, and 
\item $\mathcal{Z}_{n,{\lambda}}(q_+,q_-,q_1)\big/\mathcal{Z}_{n+1,{\lambda}}(q_+,q_-,q_1)$ is described explicitly. 
\end{itemize}
As a corollary, we get the following unchanging of the normalized generating functions under ``wall-crossing":
\begin{itemize}
\item $\mathcal{Z}_{0,\vlam}\big/\mathcal{Z}_{0,{\emptyset}}=\mathcal{Z}_{1,\vlam}\big/\mathcal{Z}_{1,{\emptyset}}=\cdots =\mathcal{Z}_{\vlam}\big/\mathcal{Z}_{{\emptyset}}$ if $\lambda^{++}=\lambda^{--}=\emptyset$, and
\item
$\Bigl(\mathcal{Z}_{0,\vlam}\big/\mathcal{Z}_{0,{\emptyset}}\Bigr)
\Big|_{q_+=q_-}=
\Bigl(\mathcal{Z}_{1,\vlam}\big/\mathcal{Z}_{1,{\emptyset}}\Bigr)
\Big|_{q_+=q_-}
=\cdots =
\Bigl(\mathcal{Z}_{\vlam}\big/\mathcal{Z}_{{\emptyset}}\Bigr)
\Big|_{q_+=q_-}$.
\end{itemize}
In this paper we will generalize the results in \cite{open_conifold} for $X_{a,b}$, where we can generalize for $X_{\Z/2\Z}$ straightforwardly. 

As in \cite{open_conifold}, the main tool is the operation so called {\it dimer shuffling}, which is nothing but the {\it mutation} in the categorical point of view. 
In \cite{3tcy}, we study dimer model descriptions for $X_{a,b}$ and mutations for them implicitly. 
On the other hand, 
\begin{itemize}
\item the dimer model for the conifold is well-known for many specialists, and 
\item the mutation for the conifold is trivial.
\end{itemize}
Once we describe the dimer model and their mutations, the argument is completely parallel to that in \cite{open_conifld}.
So what we will do mainly in this paper is describing them explicitly. 
\end{NB}

\subsection*{Acknowledgement}
The author is supported by JSPS Fellowships for Young Scientists (No.\ 19-2672).
He thanks to Hiroaki Kanno and Masahito Yamazaki for useful comments.

\subsection*{Notation}
\subsubsection{Indices}\label{subsub_001}
Let $\Zh$ denote the set of half integers and $L$ be a positive integer.
We set $I:=\Z/L\Z$ and $I_{\mathrm{h}}:=\Zh/L\Z$. The two natural projections $\Z\to I$ and $\Zh\to I_{\mathrm{h}}$ are denoted by the same symbol $\pi$.
We sometimes identify $I$ and $\Ih$ with $\{0,\ldots,L-1\}$ and $\{1/2,\ldots,L-1/2\}$ respectively.

The symbols $n$, $h$, $i$ and $j$ are used for elements in $\Z$, $\Zh$, $I$ and $\Ih$ respectively.

For $n\in\Z$ and $h\in\Zh$, we define $c(n),c(h)\in\Z$ by
\[
n=c(n)\cdot L+\pi(n),\quad
h=c(h)\cdot L+\pi(h).
\]

\subsubsection{Young diagram}\label{subsub_002}
A {\it Young diagram} $\nu$ is a map $\nu\colon\Z\to\Z$ such that $\nu(n)=|n|$ if $|n|\gg 0$ and $\nu(n)-\nu(n-1)=\pm 1$ for any $n\in\Z$. 
The map $\Zh\to \{\pm 1\}$ given by $j\mapsto\nu(j+1/2)-\nu(j-1/2)$ is also denoted by $\nu$.

By an abuse of notation, we sometimes identify $+$ and $-$ with $1$ and $-1$.

A Young diagram can be represented by a non-increasing sequence of positive integers. We fix the notation as in Figure \ref{fig15}.
\begin{figure}[htbp]
  \centering
  \input{fig15.tpc}
  \caption{$\nu=(1,1)$, ${}^{\mathrm{t}}\vspace{-2pt}\nu=(2)$.}
\label{fig15}
\end{figure}

\subsubsection{Formal variables}\label{subsub_003}
Let $q_+,q_-$ and $q_0,\ldots,q_{L-1}$ be formal variables.
We use $q_+,q_-$ and $q_1,\ldots,q_{L-1}$ for generating functions of refined invariants. 
Substituting $q_+=q_-=(q_0)^{1/2}$, we get generating functions of non-refined invariants. 

Let $P:=\Z\cdot I$ be the lattice with the basis $\{\alpha_i\mid i\in I\}$.
For an element 
$\alpha=\sum\alpha^i\cdot \alpha_i\in P\quad (\alpha^i\in\Z)$, 
we put $q^\alpha:=\prod (q_i)^{\alpha^i}$.

For $\alpha$, $\alpha'\in P$, we say $\alpha<\alpha'$ or $q^\alpha<q^{\alpha'}$ if $\alpha'-\alpha\in P^+:=\Z_{\geq 0}\cdot I$.

\section{Preliminaries}
\subsection{Affine root system}
\subsubsection{}
For $h$, $h'\in\Zh$, we define $\alpha_{[h,h']}\in P$ by
\[
\alpha_{[h,h']}:=
{\displaystyle \sum_{n=h+1/2}^{h'-1/2}\alpha_{\pi(n)}}
\]
if $h<h'$, $\alpha_{[h,h']}=1$ if $h=h'$ and $\alpha_{[h,h']}=-\alpha_{[h',h]}$ if $h>h'$.
We set
\begin{align*}
\Lambda &:=\{\alpha_{[h,h']}\in P\mid h\neq h'\},\\
\Lambda^{\mathrm{re},+}&:=\{\alpha_{[h,h']}\in \Lambda\mid h<h',\ h\not\equiv h'\ (\mathrm{mod}\,L)\}.
\end{align*}
An element in $\Lambda$ (resp. $\Lambda^{\mathrm{re},+}$) is called a {\it root} (resp. {\it positive real root}) of the affine root system of type $A_{L-1}$.

\subsubsection{}
The element $\delta:=\alpha_0+\cdots \alpha_{L-1}\in P$ is called the {\it minimal imaginary root}. 
We set
\[
\Lambda^{\mathrm{fin},+}:=\bigl\{\alpha_{[j,j']}\in \Lambda\mid 1/2\leq j<j'\leq L-1/2\bigr\}
\]
and 
\begin{equation}\label{eq_re++}
\Lambda^{\mathrm{re},+}_+:=\bigl\{\alpha_{[j,j']}+N\delta\mid \alpha_{[j,j']}\in\Lambda^{\mathrm{fin},+},\,N\geq 0\bigr\}.
\end{equation}
\begin{ex}
In the case of $L=4$, we have 
\[
\Lambda^{\mathrm{fin},+}:=\bigl\{\alpha_1,\alpha_2,\alpha_3,\alpha_2+\alpha_3,\alpha_1+\alpha_2, \alpha_1+\alpha_2+\alpha_3\bigr\}.
\]
\end{ex}

\subsubsection{}
For a root $\alpha\in \Lambda$, we take $h$ and $h'$ such that $\alpha=\alpha_{[h,h']}$ and set 
\[
j_-(\alpha):=\pi(h), \text{ and } j_+(\alpha):=\pi(h').
\]
We also put 
\[
B^\alpha:=\bigl\{(h,h')\in\bigl(\Zh\bigr)^2\,\big|\,\alpha_{[h,h']}=\alpha\bigr\}.
\]

\subsubsection{}
Let $\Theta$ denote the set of bijections $\theta\colon \Zh\to\Zh$ such that 
\begin{itemize}
\item $\theta(h+L)=\theta(h)+L$ for any $h\in\Zh$, and
\item ${\displaystyle \sum_{h=1/2}^{L-1/2}\theta(h)=\sum_{h=1/2}^{L-1/2}h}$.
\end{itemize}
\begin{ex}
In the case of $L=4$, the correspondence 
\[
\frac{1}{2}\mapsto -\frac{1}{2},\quad
\frac{3}{2}\mapsto \frac{3}{2},\quad
\frac{5}{2}\mapsto \frac{5}{2},\quad
\frac{7}{2}\mapsto \frac{9}{2}
\]
gives an elements in $\Theta$.
Let $\mu_0(\mathrm{id})$ denote this map \textup{(}see \S \ref{subsub_121} for the notation\textup{)}.
\end{ex}

\subsubsection{}
For $\theta\in \Theta$ and $i\in I$, we define $\alpha(\theta,i):=\alpha_{[\theta(n-1/2),\theta(n+1/2)]}$ ($n\in \pi^{-1}(i)$).
\begin{ex}
\[
\begin{array}{ll}
\alpha(\mathrm{id},0)=\alpha_0, \quad&
\quad\alpha(\mu_0(\mathrm{id}),0)=-\alpha_0, \\
\alpha(\mathrm{id},1)=\alpha_1, \quad&
\quad\alpha(\mu_0(\mathrm{id}),1)=\alpha_0+\alpha_1, \\
\alpha(\mathrm{id},2)=\alpha_2, \quad&
\quad\alpha(\mu_0(\mathrm{id}),2)=\alpha_2, \\ 
\alpha(\mathrm{id},3)=\alpha_3, \quad&
\quad\alpha(\mu_0(\mathrm{id}),3)=\alpha_0+\alpha_3.
\end{array}
\]
\end{ex}

\subsubsection{}
A positive real root $\alpha=\alpha_{[h,h']}$ is said to take a positive (resp. negative) value with respect to $\theta$, which is denoted by $\theta(\alpha)>0$ (resp. $<0$), if $\theta^{-1}(h)>\theta^{-1}(h')$ (resp. $\theta^{-1}(h)<\theta^{-1}(h')$). We set
\begin{equation}\label{eq_lam_theta}
\Lambda^{\mathrm{re},+}_\theta:=\Bigl\{\alpha\in \Lambda^{\mathrm{re},+}\,\Big|\,\theta(\alpha)>0\,\Bigr\}.
\end{equation}
\begin{ex}
\[
\Lambda^{\mathrm{re},+}_{\mathrm{id}}=\emptyset,\quad 
\Lambda^{\mathrm{re},+}_{\mu_0(\mathrm{id})}=\{\alpha_0\}.
\]
\end{ex}
\begin{rem}
As we mentioned in the introduction, we studied moduli spaces of representations of a noncommutative crepant resolution of $X_{L^+,L^-}$ in \cite{3tcy}. 
In this theory, the space of stability conditions can be canonically identified with $P^*\otimes \R$ and the walls are classified as follows:
\begin{itemize}
\item the walls $W_\alpha:=(\R\cdot \alpha)^{\bot}\subset P^*\otimes \R$\quad \textup{(}$\alpha\in \Lambda^{\mathrm{re},+}$\textup{)}, and 
\item the wall $W_\delta:=(\R\cdot \delta)^{\bot}$, which separates the Donaldson-Thomas and Pandharipande-Thomas chambers. 
\end{itemize}
The maps $\theta\colon \Zh\to\Zh$ as above parametrize the chambers on one side of the wall $W_\delta$.
The notation $\theta(\alpha)\gtrless 0$ respects this parametrization.
\end{rem}

\subsection{``Wall-crossing"}
\subsubsection{}\label{subsub_121}
For $i\in I$, let $\mu_i\colon \Zh\to\Zh$ be the map given by
\[
\mu_i(h)=
\begin{cases}
h-1, & \pi(h-1/2)=i,\\
h+1, & \pi(h+1/2)=i,\\
h, & \text{otherwise}.
\end{cases}
\]
For $\theta\in\Theta$, we put $\mu_i(\theta):=\theta\circ \mu_i$.
\begin{rem}
The chambers corresponding to $\theta$ and $\mu_i(\theta)$ are separated by the wall $W_{\alpha(\theta,i)}$, which is the reason for the title of this subsection. 
From the viewpoints of the affine root system, noncommutative crepant resolutions and dimer models, ``wall-crossing" corresponds to ``simple reflection", ``mutation" and ``dimer shuffling" respectively.
\end{rem}

\subsubsection{}
Let $\mathbf{i}=(i_1,i_2,\ldots)\in I^{\Z_{>0}}$ be a sequence of elements in $I$. 
For $b>0$, we define 
\[
\theta_{\mathbf{i},b}:=\mu_{i_{b-1}}(\cdots(\mu_{i_1}(\mathrm{id}))\cdots)\in\Theta,\quad
\alpha_{\mathbf{i},b}:=\alpha(\theta_{\mathbf{i},b},i_b).
\]
We say $\mathbf{i}\in I^{\Z_{>0}}$ is a {\it minimal expression} if 
$\theta_{\mathbf{i},{b}}(\alpha_{\mathbf{i},{b}})<0$
for any $b{>0}$.
For a minimal expression $\mathbf{i}$, we have
\[
\Lambda^{\mathrm{re},+}_{\theta_{\mathbf{i},b}}=\{\alpha_{\mathbf{i},1},\ldots,\alpha_{\mathbf{i},b-1}\}.
\]
\subsection{Core and quotient of a Young diagram}
\subsubsection{}\label{subsub_core}
Let $\sigma\colon I_{\mathrm{h}}\to \{\pm\}$ and $\lambda\colon\Zh\to \{\pm\}$ be maps such that $\lambda(h)=\pm \sigma (\pi (h))$ if $\pm h\gg 0$.

We define integers $c_{\lambda}[j]$ and Young diagrams $\lambda^{[j]}$ for $j\in I_{\mathrm{h}}$ by 
\[
\lambda(h)=\lambda^{[\pi(h)]}\bigl(\sigma(j(h))\cdot (c(h)-c_{\lambda}[\pi(h)]+1/2)\bigr).
\]
\begin{rem}
In the case $\sigma\equiv +$ and $\sum c_{\lambda}[j]=0$, the sequence $\bigl(c_{\lambda}[j]\bigr)$ of integers and the sequence $\bigl(\lambda^{[j]}\bigr)$ of Young diagrams are called the \textup{$L$-core} and the \textup{$L$-quotient} of the Young diagram $\lambda$.
\end{rem}

\subsubsection{}
We put
\begin{equation}\label{eq_Bpm1}
B^{\alpha,\pm}_{\sigma,\lambda}:=\{(h,h')\in B^\alpha\mid -\lambda(h)\sigma(h)=\lambda(h')\sigma(h')=\pm\}.
\end{equation}
\begin{lem}\label{lem_14}
\[
\bigl|B^{\alpha,+}_{\sigma,\lambda}\bigr|-\bigl|B^{\alpha,-}_{\sigma,\lambda}\bigr|=\alpha^0+c_{\lambda}[j_-(\alpha)]-c_{\lambda}[j_+(\alpha)].
\]
\end{lem}
\begin{proof}
We write simply $j_\pm$ for $j_\pm(\alpha)$.
Note that we have
\[
B^\alpha=\Bigl\{\bigl(cL+j_-,(c+\alpha^0)L+j_+\bigr)\,\Big|\, c\in\Z\Bigr\}. 
\]
For an integer $N$, we set 
\[
B^\alpha_N:=\Bigl\{\bigl(cL+j_-,(c+\alpha^0)L+j_+\bigr)\,\Big|\,
c\in[-N,N-1]\Bigr\}.
\]
Take a sufficiently large integer $N$, then we have
\[
B^{\alpha,+}_{\sigma,\lambda},B^{\alpha,-}_{\sigma,\lambda}\subset B^\alpha_N
\]
and so
\begin{align*}
&\bigl|B^{\alpha,+}_{\sigma,\lambda}\bigr|-\bigl|B^{\alpha,-}_{\sigma,\lambda}\bigr|\\
&=
-\sharp\bigl\{(h,h')\in B^\alpha_N\,\big|\,\lambda(h)\sigma(h)=+\bigr\}
+\sharp\bigl\{(h,h')\in B^\alpha_N\,\big|\,\lambda(h')\sigma(h')=+\bigr\}\\
&=-\sharp\bigl\{c\in[-N,N-1]\,\big|\, \lambda^{[j_-]}\bigl(\sigma(j_-)\cdot(c-c_\lambda[j_-]+1/2)\bigr)=\sigma(j_-)\bigr\}\\
&\ \quad +\sharp\bigl\{c\in[-N,N-1]\,\big|\, \lambda^{[j_+]}\bigl(\sigma(j_+)\cdot(c+\alpha^0-c_\lambda[j_+]+1/2)\bigr)=\sigma(j_+)\bigr\}\\
&=-(N-c_\lambda[j_-]-1/2)+(N+\alpha^0-c_\lambda[j_+]-1/2)\\
&=\alpha^0+c_{\lambda}[j_-]-c_{\lambda}[j_+].
\end{align*}
\end{proof}
For $\sigma, \lambda, \theta$ and $i$, we put
\begin{equation}\label{eq_Bpm2}
B^{i,\pm}_{\plt}:=\bigl\{n\in \pi^{-1}(i)\mid (\theta(n-1/2),\theta(n+1/2))\in B^{\alpha(\theta,i),\pm}_{\sigma,\lambda}\bigr\}.
\end{equation}

\section{Dimer model}\label{sec_dimer_model}
\subsection{Dimer configurations}
\subsubsection{}
We fix the following data:
\begin{itemize}
\item a map $\sigma\colon I_{\mathrm{h}}\to \{\pm\}$, 
\item a map $\lambda\colon\Zh\to \{\pm\}$ such that $\lambda(h)=\pm \sigma (\pi (h))$ for $\pm h\gg 0$,
\item a pair of Young diagrams $\nu=(\nu_+,\nu_-)$,
\item a bijection $\theta\colon \Zh\to\Zh$ in $\Theta$.
\end{itemize}
We put $\tsig:=\sigma\circ\pi\circ\theta$, $\tlam:=\lambda\circ\theta$ and $L_\pm:=|\sigma^{-1}(\pm)|$.

\subsubsection{}\label{subsec_graph}
We consider the following graph in the $(x,y)$-plane. 

First, we set 
\begin{align}
H(\sigma,\theta)&:=\bigr\{n\in\Z\mid \tsig(n-1/2)=\tsig(n+1/2)\bigl\},\quad I_H(\sigma,\theta):=\pi(H(\sigma,\theta)),\label{eq_IH}\\
S(\sigma,\theta)&:=\bigr\{n\in\Z\mid \tsig(n-1/2)\neq \tsig(n+1/2)\bigl\},\quad I_S(\sigma,\theta):=\pi(S(\sigma,\theta))\label{eq_IS}
\end{align}
and for $n\in H(\sigma,\theta)$ we put $\tsig(n):=\tsig(n\pm1/2)$.

The set of the vertices is given by 
\begin{align*}
\mathcal{V}:=&\, 
\bigr\{(n,m)\mid n\in S(\sigma,\theta),\ \text{$n-m$: odd}\bigl\}\\
&\sqcup\{(n-1/2,m)\mid n\in H(\sigma,\theta),\ \text{$n-m$: odd}\}\\
&\sqcup\{(n+1/2,m)\mid n\in H(\sigma,\theta),\ \text{$n-m$: odd}\},
\end{align*}
which are denoted by $v(n,m)$, $\vl(n-1/2,m)$ and $\vr(n+1/2,m)$ respectively.

The set of the edges is given by 
\[
\mathcal{E}:= 
\bigl\{\eh(n,m)\mid n\in H(\sigma,\theta),\ \text{$n-m$\,:\,odd}\bigl\}\,\sqcup\,\bigr\{\es(h,k)\mid h,k\in \Zh\bigr\},
\]
where 
\begin{itemize}
\item $\eh(n,m)$ connects $\vl(n-1/2,m)$ and $\vr(n+1/2,m)$,
\item $\es(h,k)$ connects $v(h-1/2,k+1/2)$ or $\vr(h,k+1/2)$ and $v(h+1/2,k-1/2)$ or $\vl(h,k-1/2)$ if $h-k$ is even, and 
\item $\es(h,k)$ connects $v(h-1/2,k-1/2)$ or $\vr(h,k-1/2)$ and $v(h+1/2,k+1/2)$ or $\vl(h,k+1/2)$ if $h-k$ is odd.
\end{itemize}
We put
\begin{equation}\label{eq_F}
\mathcal{F}:=\{(n,m)\in\Z^2\mid n+m:\text{even}\},\ \mathcal{F}_i:=\{(n,m)\in\mathcal{F}\mid n\in\pi^{-1}(i)\}
\end{equation}
for $i\in I$. 
Note that $\mathcal{E}$ divides the plain into disjoint hexagons and quadrilaterals, where the hexagons are parametrized by the set
\[
\mathcal{F}_{\mathrm{H}}:=\{(n,m)\in\mathcal{F}\mid n\in H(\sigma,\theta)\}
\]
and the quadrilaterals are parametrized by the set
\[
\mathcal{F}_{\mathrm{S}}:=\{(n,m)\in\mathcal{F}\mid n\in S(\sigma,\theta)\}.
\]
For $(n,m)\in\mathcal{F}$, let $f(n,m)$ denote the corresponding hexagon or quadrilateral.
\begin{ex}\label{ex1}
In Figure \ref{fig1}, we show the graph associated with $L=3$, $\sigma$ given by 
\[
\sigma(1/2)=+.\quad
\sigma(3/2)=-.\quad
\sigma(5/2)=-
\]
and $\theta=\mathrm{id}$ \textup{(}$L_+=1, L_-=2$\textup{)}.
\begin{figure}[htbp]
  \centering
  \input{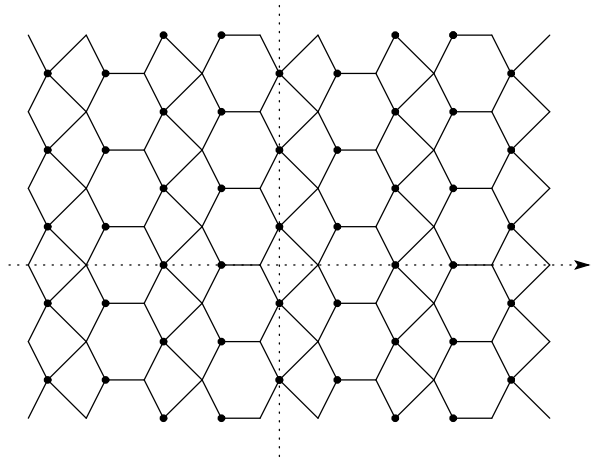}
  \caption{The graph and $\mathcal{V}_+$.}
\label{fig1}
\end{figure}
\end{ex}

\subsubsection{}\label{subsec_bipartite}
We set 
\begin{align*}
\mathcal{V}_\pm:=&\, 
\bigr\{v(n,m)\mid \tsig(n+1/2)=\pm\bigl\}\\
&\sqcup\{\vl(n-1/2,m)\mid \tsig(n)=\mp\}
\sqcup\{\vr(n+1/2,m)\mid \tsig(n)=\pm\}.
\end{align*}
Note that $\mathcal{V}=\mathcal{V}_+\sqcup\mathcal{V}_-$ and each element in $\mathcal{E}$ connects an element in $\mathcal{V}_+$ and an element in $\mathcal{V}_-$ (see Figure \ref{fig1} for example). 

A {\it perfect matching} is a subset of $\mathcal{E}$ giving a bijection between $\mathcal{V}_+$ and $\mathcal{V}_-$.

\subsubsection{}
We define the map $F_{\plt}\colon\Z\to\Z$ by $F_{\plt}(0)=0$ and
\begin{equation}\label{eq_211}
F_{\plt}(n)=F_{\plt}(n-1)-\tlam(n-1/2).
\end{equation}
For $k\in\Zh$, we set
\begin{align*}
\mathcal{P}_{\plt}^{k,\pm}:=&\,\Bigr\{\eh\bigr(n,F_{\plt}(n)+2k\bigl)\,\Big|\, n\in\Z,\ \tsig(n)=\mp\Bigl\}\\
&\sqcup \Bigr\{\es\bigr(h,(F_{\plt}(h-1/2)+F_{\plt}(h+1/2))/2+2k\bigl)\,\Big|\,h\in\Zh,\ \tsig(h)=\pm\Bigl\}.
\end{align*}
For a Young diagram $\eta$, let $\mathcal{P}_{\plt}^{\eta}$ denote the following perfect matching:
\[
\mathcal{P}_{\plt}^{\eta}:=\bigsqcup_{k\in\Zh}\mathcal{P}_{\plt}^{k,{\eta}(k)}.
\]  
\begin{ex}\label{ex2}
In Figure \ref{fig4}, we show the perfect matching associated with $\sigma$ as is Example \ref{ex1}, $\theta=\mathrm{id}$, $\eta=\emptyset$ and $\lambda$ given by
\[
\lambda(h)=
\begin{cases}
+, & h=-5/2,\\
-, & h=1/2,\\
\mathrm{sgn}(h)\sigma(h), & \textup{otherwise}.
\end{cases}
\]
\begin{figure}[htbp]
  \centering
  \input{fig4.tpc}
  \caption{$\{f(n,F_{\sigma,\lambda,\mathrm{id}}(n))\mid n\in\Z\}$ and $\mathcal{P}_{\sigma,\lambda,\mathrm{id}}^{\emptyset}$}
\label{fig4}
\end{figure}
\end{ex}

\begin{NB}
\subsubsection{}
We set 
\begin{align*}
W(\sigma,\theta)&:=\{j\in\Zh\mid p\circ\theta(1/2)=p\circ\theta(j)\},\\
B(\sigma,\theta)&:=\{j\in\Zh\mid p\circ\theta(1/2)\neq p\circ\theta(j)\},
\end{align*}
and 
\begin{align*}
\tE_{\plt}^{k,+}:=&\,\Bigr\{\eh\bigr(n,F_{\plt}(n)+2k\bigl)\,\Big|\, n\in\Z,\ n+1/2\in W(p.\theta)\Bigl\}\\
&\sqcup \Bigr\{\es\bigr(j,(F_{\plt}(j-1/2)+F_{\plt}(j+1/2))/2+2k\bigl)\,\Big|\,j\in\Zh,\ i\in B(p.\theta)\Bigl\},\\
\tE_{\plt}^{k,+}:=&\,\Bigr\{\eh\bigr(n,F_{\plt}(n)+2k\bigl)\,\Big|\, n\in\Z,\ n+1/2\in B(p.\theta)\Bigl\}\\
&\sqcup \Bigr\{\es\bigr(j,(F_{\plt}(j-1/2)+F_{\plt}(j+1/2))/2+2k\bigl)\,\Big|\,j\in\Zh,\ i\in W(p.\theta)\Bigl\},
\end{align*}
and
\[
\tE_{\plt}^{\nu_\pm}:=\bigsqcup_{k\in\Zh}\tE_{\plt}^{k,{\nu_\pm}(k)}.
\]  
Note that for any vertex $v\in \mathcal{V}$ we have the unique edge $e\in \tE_{\plt}^\mu$ with $v$ as its end point. 
A subset of $\mathcal{E}$ satisfying this condition is said to be a {\it perfect matching} of the graph. 
\end{NB}

\subsubsection{}
Let $\mathcal{P}^\pm_{\plt}$ be the following perfect matching:
\[
\mathcal{P}^\pm_{\plt}:=
\bigl\{\eh(n,m)\mid \tsig(n)=\mp\bigr\}
\sqcup
\bigl\{\es(h,k)\mid \tsig(h)=\pm,\ h\cdot\tlam(h)-k:\text{even}\bigr\}.
\]

\subsubsection{}
\begin{defn}
A perfect matching $D$ is said to be {\it a dimer configuration of type $(\sigma,\lambda,\nu,\theta)$} if $D$ coincides with 
$\mathcal{P}_{\plt}^{\nu_\pm}$ in the area $\{\pm x>m\}$ and
$\mathcal{P}^\pm_{\plt}$ in the area $\{\pm y>m\}$
 for $m\gg 0$.
\end{defn}
\begin{rem}\label{rem_brane}
A dimer configuration of type $(\sigma,\vec{\emptyset},\vec{\emptyset},\mr{id})$ is ``a perfect matching congruent to the canonical perfect matching'' in the terminology of \cite{ncdt-brane}.
\end{rem}

\subsubsection{}\label{subsub_217}
For $f\in \mathcal{F}$, let $\partial f\subset \mathcal{E}$ denote the set of edges surrounding the face $f$. 
By moving around the face $f$ clockwisely, we can determine an orientation for each element in $\partial f$.
Let $\partial^{\pm} f\subset \partial f$ denote the subset of edges starting from elements in $\mathcal{V}_\pm$. 

For an edge $e\in\mathcal{E}$, let $f^\pm(e)$ denote the unique face such that $e\in \partial^{\pm} f^\pm(e)$.

\subsection{Weights}
\subsubsection{}\label{subsub_221}
For $h\in\Zh$, we define the monomials $w_{\sigma,\lambda}(h)$ by the following conditions:
\[
w_{\sigma,\lambda}(h)=
\begin{cases}
\bigl(Q_{\sigma(h)}\bigr)^{c(h)-c_\lambda[j(h)]}q_{\sigma(h)}^{(j(h))} & h\gg 0.\\
\bigl(Q_{-\sigma(h)}\bigr)^{c(h)-c_\lambda[j(h)]}q_{-\sigma(h)}^{(j(h))} & h\ll 0,
\end{cases}
\]
and 
\[
w_{\sigma,\lambda}(h)/w_{\sigma,\lambda}(h-L)=q_{\lambda(h)}\cdot q_{\lambda(h-L)}\cdot q_1\cdot\cdots\cdot q_{L-1},
\]
where 
\[
Q_\pm:=(q_\pm)^2\cdot q_1\cdot\cdots\cdot q_{L-1},\ 
q_\pm^{(j)}:=q_\pm\cdot  q_1\cdot\cdots\cdot q_{j-1/2}.
\]
Note that for $h\neq h'$ we have
\begin{equation}\label{eq_221}
w_\lambda(h')/w_\lambda(h)\Big|_{q_+=q_-=(q_0)^{1/2}}=q^{\alpha_{[h,h']}}.
\end{equation}
\begin{ex}\label{ex3}
In Figure \ref{fig16}, we show the weight $w_{\sigma,\lambda}$ for $\sigma$ and $\lambda$ as in Example \ref{ex2}.
\begin{figure}[htbp]
  \centering
  \input{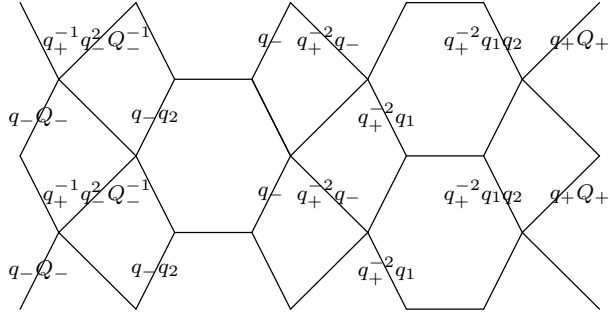}
  \caption{the weight $w_{\sigma,\lambda}$}
\label{fig16}
\end{figure}
\end{ex}

\subsubsection{}
For an edge $e\in\mathcal{E}$, we associate the weight $w_{\plt}(e)$ by 
\begin{align}
w_{\plt}(\es(h,k))&:=
\begin{cases}
w_{\sigma,\lambda}(\theta(h))^{\tsig(h)\cdot \tlam(h)}, & h\cdot\tlam(h)-k:\text{odd}, \notag \\
1, & h\cdot\tlam(h)-k:\text{even},
\end{cases}\\
w_{\lambda,\sigma,\theta}(\eh(n,m))&:=1.\label{eq_weight}
\end{align}

\subsubsection{}
Fix $\sigma$ and $\lambda$, then the set
\[
\bigsqcup_{\alpha\in\Lambda^{\mathrm{re},+}}B^{\alpha,-}_{\sigma,\lambda}
\]
is finite. 
We set
\begin{align}
F^\alpha_{\sigma,\lambda}&:=\prod_{(h,h')\in B^{\alpha,-}_{\sigma,\lambda}}w_{\sigma,\lambda}(h')/w_{\sigma,\lambda}(h),\label{eq_F1}\\
F^\theta_{\sigma,\lambda}&:=\prod_{
\begin{subarray}{c}
\alpha\in \Lambda^{\mathrm{re},+};\ 
\theta(\alpha)<0,\\
\sigma(j^-(\alpha))\neq\sigma(j^+(\alpha)).
\end{subarray}}
F^\alpha_{\sigma,\lambda}.\label{eq_222}
\end{align}

\subsubsection{}
Note that for a dimer configuration $D$ of type $(\plmt)$ we have only a finite number of $e\in D$ such that $w_{\plt}(e)\neq 1$.
\begin{defn}
For a dimer configuration $D$ of type $(\plmt)$, we define the weight $w_{\plt}(D)$ by 
\begin{equation}\label{eq_def_wt}
w_{\plt}(D):=F^\theta_{\sigma,\lambda}\cdot \prod_{e\in D}w_{\plt}(e).
\end{equation}
\textup{(}See \eqref{eq_weight} and \eqref{eq_222} for the notations.\textup{)}
\end{defn}
\begin{rem}
We will define the generating function $\mathcal{Z}_{\plmt}$ by the sum of weighs of all dimer configurations of type $(\plmt)$\footnote{We will leave the definition of the generating function until \S \ref{subsec_gen_func} since we will use Proposition \ref{prop_24} to prove that the number of dimer configurations of type $(\plmt)$ is finite.}. 
\end{rem}

\subsubsection{}
For a finite subset $\mathcal{E}'\subset\mathcal{E}$, we put 
\[
w_{\plt}(\mathcal{E}'):=\prod_{e\in\mathcal{E}'}w_{\plt}(e)
\]
and for a face $f\in\mathcal{F}$ we put 
\begin{equation}\label{eq_223}
w_{\plt}(f):=w_{\plt}(\partial^- f)/w_{\plt}(\partial^+ f).
\end{equation}
For an integer $n$ we set 
\[
w_{\plt}(n):=w_{\sigma,\lambda}(\theta(n+1/2))/w_{\sigma,\lambda}(\theta(n-1/2)),
\]
then we have
\[
w_{\plt}(f(n,m))=w_{\plt}(n)
\]
for any $(n,m)\in\mathcal{F}$.
By \eqref{eq_221}, we have
\[
w_{\plt}(n)\big|_{q_+=q_-=(q_0)^{1/2}}=q^{\alpha(\theta,i)}.
\]

\section{Viewpoint from noncommutative crepant resolution}\label{sec_crystal}
\subsection{Noncommutative crepant resolution}\label{subsec_nccr}
Let $\Gamma$ be a lattice in the $(x,y)$-plane generated by $(L,0)$ and $(0,2)$. The graph given in \S \ref{subsec_graph} is invariant under the action of $\Gamma$ and so gives a graph on the torus $\R^2/\Gamma$. 
This gives a quiver with a potential $A=(Q_{\sigma,\theta},w_{\sigma,\theta})$ as in \cite{3tcy}.
The vertices of $Q_{\sigma,\theta}$ are parametrized by $I$ and the arrows are given by 
\[
\Biggr(\,\bigsqcup_{j\in\Ih}h_j^+\Biggr)
\,\sqcup\,
\Biggr(\,\bigsqcup_{j\in\Ih}h_j^-\Biggr)
\,\sqcup\,
\Biggr(\bigsqcup_{i\in I_H(\sigma,\theta)}r_i\Biggr)
\]
(see \eqref{eq_IH} for the notation).
Here $h^+_j$ (resp. $h^-_j$) is an edge from $j-1/2$ to $j+1/2$ (resp. from $j+1/2$ to $j-1/2$) and $r_i$ is an edge from $i$ to itself.
See \cite[\S 1.2]{3tcy} for the definition of the potential $w_{\sigma,\theta}$.\begin{ex}
In Figure \ref{fig27}, we show the quiver $Q_{\sigma,\mathrm{id}}$ for $\sigma$  as in Example \ref{ex2}.
\begin{figure}[htbp]
  \centering
  \input{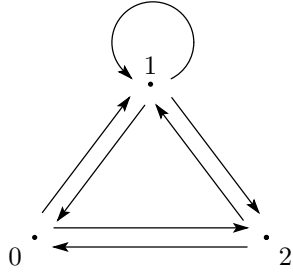}
  \caption{the quiver $w_{\sigma,\theta}$}
\label{fig27}
\end{figure}
\end{ex}
\begin{rem}
\begin{itemize}
\item The center of $A$ is isomoriphic to $R:=\C[\mr{x},\mr{y},\mr{z},\mr{w}]/(\mr{x}\mr{y}=\mr{z}^{L_+}\mr{w}^{L_-})$.
In \cite[Theorem 1.14 and 1.19]{3tcy}, we showed that $A$ is a noncommutative crepant resolution of $X=\mathrm{Spec}\, R$.
\item The affine $3$-fold $X$ is toric. 
In fact, 
\[
T=\mathrm{Spec}\,\tilde{R}:=\mathrm{Spec}\,\C[\mr{x}^\pm,\mr{y}^\pm,\mr{z}^\pm,\mr{w}^\pm]/(\mr{x}\mr{y}=\mr{z}^{L_+}\mr{w}^{L_-})\subset X
\]
is a $3$-dimensional torus.
\end{itemize}
\end{rem}

\subsection{Dimer model and noncommutative crepant resolution}\label{subsec_dimer_nccr}

\subsubsection{}
We will construct an $A$-module $M(D)$ for a dimer configuration $D$. 
Let $V_i=V_i(D)$ ($i\in I$) be vector space with the following basis: 
\[
\bigl\{b[D;x,y,z]\mid (x,y)\in\mathcal{F}_i,\ z\in\Z_{\geq 0}\bigr\}
\]
(see \eqref{eq_F} for the notation).
We define the map $h_j^\pm\colon V_{j\mp1/2}\to V_{j\pm1/2}$ by 
\[
h_j^\pm(b[D;x,y,z])=\begin{cases}
b[D;x\pm1,y-\tsig(j),z] & \es(x\pm1/2,y-\tsig(j)/2)\notin D,\\
b[D;x\pm1,y-\tsig(j),z+1] & \es(x\pm1/2,y-\tsig(j)/2)\in D,
\end{cases}
\]
and $r_i\colon V_i\to V_i$ by 
\[
r_i(b[D;x,y,z])=
\begin{cases}
b[D;x,y+\tsig(j),z] & \eh(x,y+\tsig(j)/2)\notin D,\\
b[D;x,y+\tsig(j),z+1] & \eh(x,y+\tsig(j)/2)\in D.
\end{cases}
\]
\begin{figure}[htbp]
  \centering
  \input{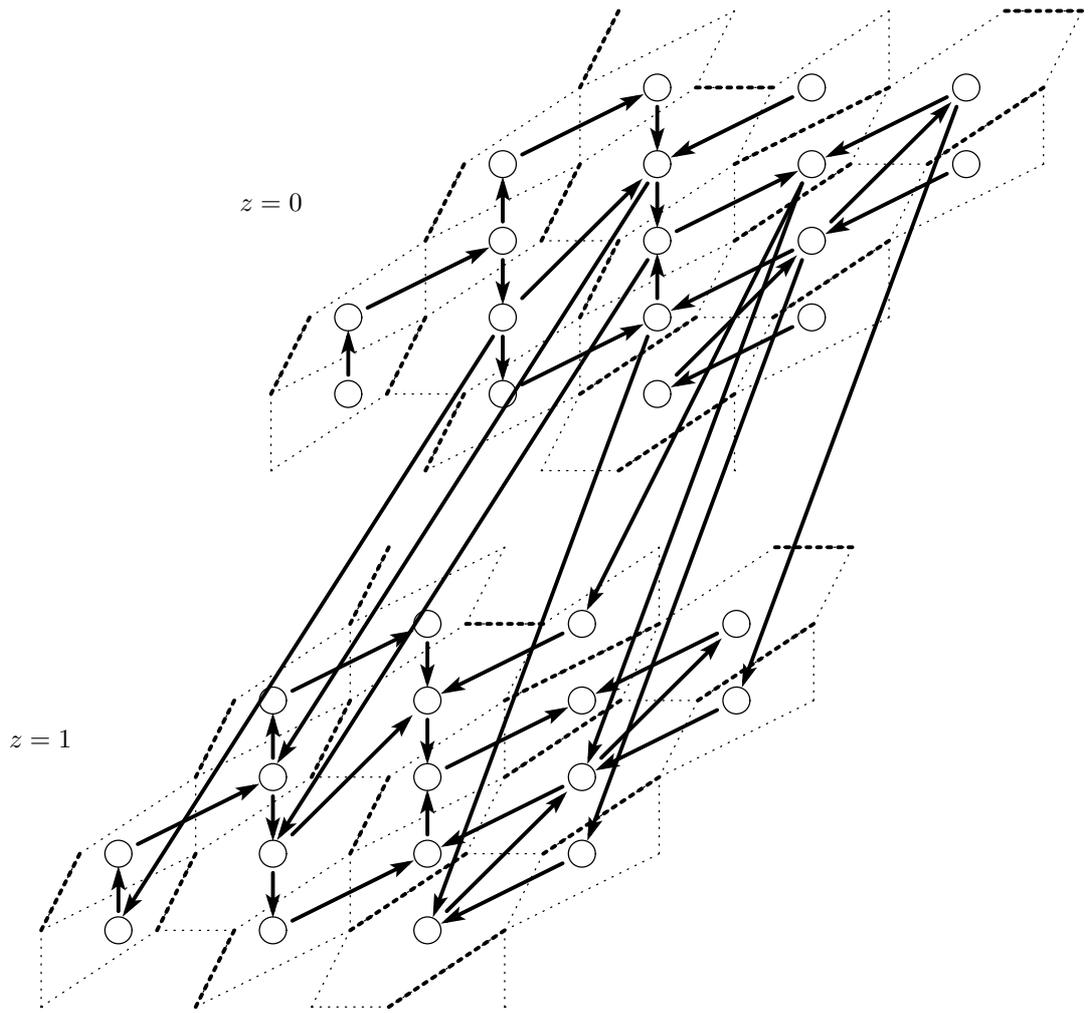}
  \caption{an example of $M(D)$}
\label{}
\end{figure}

\begin{NB}
\subsubsection{}\label{subsub_222}
A subset $\Gamma$ of $\mathcal{E}$ is said to be an {\it alternative cycle} with respect to a dimer configuration $D$ if 
\begin{itemize}
\item $\Gamma$  gives a closed zigzag curve without self-intersection, and 
\item along the closed zigzag curve, the elements in $\Gamma\cap D$ and the elements in $\Gamma -D$ appear alternatively.
\end{itemize}
\end{NB}

\subsubsection{}\label{subsub_322}
Let $\mathcal{C}\subset\mathcal{E}$ be a subset which gives a closed zigzag curve without self-intersection.
By moving along the zigzag curve clockwisely, we can determine an orientation for each element in $\mathcal{C}$.
Let $\mathcal{C}^{\pm}\subset\mathcal{C}$ denote the subset of edges starting from elements in $\mathcal{V}_\pm$.

Let $D$ be a dimer configuration of type $(\plmt)$.
A subset $\mathcal{C}$ as above is said to be a {\it positive cycle} (resp. {\it negative cycle}) with respect to $D$ if $\mathcal{C}\cap D=\mathcal{C}^+$ (resp. $\mathcal{C}^-$).

\subsubsection{}\label{subsub_224}
Given a dimer configuration $D$ and a positive cycle $\mathcal{C}$ with respect to $D$, let $D_{\mathcal{C}}$ be the dimer configuration given by 
\[
D_{\mathcal{C}}=(D\backslash {\mathcal{C}}^+)\cup {\mathcal{C}}^-.
\]
Then we can check the following lemma:
\begin{lem}\label{lem_alternative_cycle}
The surjection $M(D)\to M({D_{\mathcal{C}}})$ given by 
\[
b[D;x,y,z]\mapsto\begin{cases}
0 & (x,y)\in {\mathcal{C}}^\circ,\ z=0\\
b[D_{\mathcal{C}};x,y,z-1] & (x,y)\in {\mathcal{C}}^\circ,\ z\geq 1\\
b[D_{\mathcal{C}};x,y,z] & (x,y)\notin {\mathcal{C}}^\circ
\end{cases}
\]
is a homomorphism of $A$-modules, where ${\mathcal{C}}^\circ$ is the interior of the closed zigzag curve. Moreover we have 
\[
w_{\plt}(D_{\mathcal{C}})=w_{\plt}(D)\cdot \prod_{f\in{\mathcal{C}}^\circ}w_{\plt}(f).
\]
\end{lem}

\subsection{Crystal melting interpretation}\label{subsec_crystal}
In this subsection, we show that a dimer configuration of type $(\plmt)$ corresponds to a (torus invariant) quotient $A$-module of the $A$-module $M^{\mathrm{max}}=M^{\mathrm{max}}_{\plmt}$.
In physicists' terminology, studying such quotient modules is called {\it the crystal melting model} (see \cite{ooguri-yamazaki}) and $M^{\mathrm{max}}$ is called the {\it grand state} of the model.

\subsubsection{}
We define a Young diagram $G_{\sigma,\lambda,\theta}\colon\Z\to\Z$ by the following conditions:
\begin{itemize}
\item
$G_{\sigma,\lambda,\theta}(n)=|n|$ if $|n|\gg 0$, and
\item 
$G_{\sigma,\lambda,\theta}(n)=G_{\sigma,\lambda,\theta}(n-1)+\tsig(n-1/2)\tlam(n-1/2)$ for any $n$.
\end{itemize}
We define a map $G_{\plt}\colon \mathcal{F}\to\Z$ by
\begin{equation}\label{eq_G}
G_{\plt}(n,m):=G(n)_{\sigma,\lambda,\theta}+2\cdot |m-F_{\plt}(n)|,
\end{equation}
where $F_{\plt}(n)$ is given in \eqref{eq_211}.

\begin{ex}
In the case as Example \ref{ex2}, we have
\[
(G_{\sigma,\lambda,\mathrm{id}}(n))_{n\in\Z}=(\ldots,6,5,4,3,4,3,2,1,2,3,4,5,6,\ldots)
\]
and $G_{\sigma,\lambda,\mathrm{id}}(n,m)$ is given in Figure \ref{fig5}.
\begin{figure}[htbp]
  \centering
  \input{fig5.tpc}
  \caption{$G_{\sigma,\lambda,\mathrm{id}}(n,m)$}
\label{fig5}
\end{figure}
\end{ex}

\subsubsection{}
We define two maps $F_{\plt}^\pm\colon\Z\to\Z$ by the following conditions:
\begin{itemize}
\item
$F_{\plt}^\pm(n)=F_{\plt}(n)$ if $\pm n\gg 0$, and
\item $F_{\plt}^\pm(n)=F_{\plt}^\pm(n-1)\mp \tsig(n-1/2)$ for any $n$.
\end{itemize}
Then we define two maps $G_{\plt}^{\nu_\pm,\pm}\colon \mathcal{F}\to\Z$ by
\begin{equation}\label{eq_Gpm}
G_{\plt}^{\nu_\pm,\pm}(n,m):=\nu_{\pm}\bigl(m-F_{\plt}^\pm(n)\bigr)\pm n.
\end{equation}
\begin{ex}
In Figure \ref{fig9_11}, we show $G_{\sigma,\lambda,\mathrm{id}}^{\emptyset,+}$ and $G_{\sigma,\lambda,\mathrm{id}}^{\square,-}$ for $\sigma$ and $\lambda$ as in Example \ref{ex2}.
\begin{figure}[htbp]
  \centering
  \input{fig9.tpc}
  \input{fig11.tpc}
  \caption{$G_{\plt}^{\emptyset,+}$ and $G_{\plt}^{\square,-}$}
\label{fig9_11}
\end{figure}
\end{ex}

\subsubsection{}
We define a map $G_{\plt}^{\nu}\colon \mathcal{F}\to\Z$
by 
\[
G_{\plt}^{\nu}(n,m):=\max\Bigl(G_{\plt}(n,m),G_{\plt}^{\nu_+,+}(n,m),G_{\plt}^{\mu_-,-}(n,m) \Bigr).
\]

We can verify that 
\[
G_{\plt}^{\nu}(f^+(e))=G_{\plt}^{\nu}(f^-(e))+1\ \text{or}\ G_{\plt}^{\nu}(f^-(e))-3.
\]
for an edge $e\in\mathcal{E}$ (see \S \ref{subsub_217} for the notations).  
We define a perfect matching $D^{\mathrm{max}}=D^{\mathrm{max}}_{\plmt}$ by 
\[
e\in D^{\mathrm{max}}\iff G_{\plt}^{\nu}(f^+(e))=G_{\plt}^{\nu}(f^-(e))-3.
\]
Let $M^{\mathrm{max}}=M^{\mathrm{max}}_{\plmt}:=M(D^{\mathrm{max}}_{\plmt})$ denote the corresponding $A$-module.
\begin{ex}
In Figure \ref{fig14}, we show $G_{\sigma,\lambda,\mathrm{id}}^{\vec{\emptyset}}$ and $D^{\mathrm{max}}_{\sigma,\lambda,\vec{\emptyset},\mathrm{id}}$ for $\sigma$ and $\lambda$ as in Example \ref{ex2}.
\begin{figure}[htbp]
  \centering
  \input{fig14.tpc}
  \caption{$G_{\sigma,\lambda,\mathrm{id}}^{\vec{\emptyset}}$ and $D^{\mathrm{max}}_{\sigma,\lambda,\vec{\emptyset},\mathrm{id}}$}
\label{fig14}
\end{figure}
\end{ex}
\begin{rem}
The graph of the map $m\mapsto G_{\plt}^{\nu}(n,m)$ determines a Young diagram. This is what we denote by $\mathcal{V}_{\mathrm{min}}(n)$ in \textup{\cite[\S 3.1]{NCDTviaVO}}.
\end{rem}

\subsubsection{}
\begin{lem}\label{lem_no_cycle}
There is no positive cycle with respect to $D^{\mathrm{max}}$.
\end{lem}
\begin{proof}
Assume that we have a positive cycle $\mathcal{C}$. 
For an edge $e\in \partial \mathcal{C}$, let $f_{\mathrm{in}}(e)$ (resp. $f_{\mathrm{out}}(e)$) be the unique face such that $e\in \partial f_{\mathrm{in}}(e)$ and $f_{\mathrm{in}}(e)\in \mathcal{C}^\circ$ (resp. $e\in \partial f_{\mathrm{out}}(e)$ and $f_{\mathrm{out}}(e)\notin \mathcal{C}^{\mathrm{\circ}}$).
Then we have
\begin{equation}\label{eq_a-cycle}
G_{\plt}^{\nu}(f_{\mathrm{in}}(e))>
G_{\plt}^{\nu}(f_{\mathrm{out}}(e)).
\end{equation}
Take a face $(n,m)\in \mathcal{C}^\circ$.
If $G_{\plt}^{\nu}(n,m)=G_{\plt}^{\nu,\pm}(n,m)$, then 
\[
(n\pm n',F^\pm_{\plt}(n\pm n')-F^\pm_{\plt}(n)+m)\in \mathcal{C}^\circ
\]
for any $n'\geq 0$ by \eqref{eq_Gpm} and \eqref{eq_a-cycle}, and this is a contradiction. 
On the other hands, if $G_{\plt}^{\mu}(n,m)=G_{\plt}(n,m)$ and $\pm m\mp F_{\plt}(n)\geq 0$, then $(n,m\pm m')\in \mathcal{C}^\circ$ for any $m'\geq 0$ by \eqref{eq_G} and \eqref{eq_a-cycle}, and this is also a contradiction. 
Hence the claim follows.
\end{proof}

\subsubsection{}
For a map $H\colon\mathcal{F}\to \Z_{\geq 0}$, let $V^H_i\subset V_i(D^{\mathrm{max}})$ ($i\in I$) be the subspace spanned by the following elements:
\[
\bigl\{b[D^{\mathrm{max}};x,y,z]\mid (x,y)\in \mathcal{F}_i,\ z\geq H(x,y)\bigr\}.
\]
The following proposition gives a one-to-one correspondence between dimer configurations of type $(\plmt)$ and finite dimensional quotient modules of $M^{\mathrm{max}}_{\plmt}$.
\begin{prop}\label{prop_24}
Given a monomial $\mathbf{q}$, we have a natural bijection between the following sets:
\begin{itemize}
\item the set of dimer configurations of type $(\plmt)$ with weight $\mathbf{q}$, and 
\item the set of maps $H\colon\mathcal{F}\to \Z_{\geq 0}$ satisfying the following conditions:
\begin{itemize}
\item $H(f)=0$ except for only a finite number of $f\in\mathcal{F}$,
\item $(V_i^H)_{i\in I}$ is stable under the action of $A$, and
\item $w_{\plt}(D^{\mathrm{max}})\cdot \displaystyle\prod_{f}w_{\plt}(f)^{H(f)}=\mathbf{q}$.
\end{itemize}
\end{itemize}
\end{prop}
\begin{proof}
Let $D$ be a dimer configuration of type $(\plmt)$.
By Lemma \ref{lem_no_cycle}, $(D\cup D^{\mathrm{max}})\backslash (D\cap D^{\mathrm{max}})$ is a disjoint union $\sqcup\,\mathcal{C}_\gamma$ of a finite number of positive cycles. 
We define a map $H_D\colon \mathcal{F}\to \Z_{\geq 0}$ by 
\[
H_D(f):=\sharp\{\mathcal{C}_\gamma\mid f\in\mathcal{C}_\gamma^{\circ}\}.
\]
Then we can verify the claim using Lemma \ref{lem_alternative_cycle}.
\end{proof}
\begin{rem}
The graph of the map $m\mapsto G_{\plt}^{\nu}(n,m)+2H(n,m)$ determines a Young diagram. This is what we denote by $\mathcal{V}(n)$ in \textup{\cite[\S 3.1]{NCDTviaVO}}.
\end{rem}

\subsection{Generating function}\label{subsec_gen_func}
From the description given by Proposition \ref{prop_24}, we can verify that, fixing a monomial $\mathbf{q}$, we have only a finite number of dimer configurations of type $(\plmt)$ with weight $\mathbf{q}$.
\begin{defn}\label{defn_gen_func}
We define the generating function by 
\[
\mathcal{Z}_{\plmt}=\mathcal{Z}_{\plmt}\bigl(\vec{q}\bigr):=\sum_{D}w_{\plt}(D),
\]
where the sum is taken over all dimer configurations of type $(\plmt)$.
In particular, we put 
\[
\mathcal{Z}_{\sigma,\lambda,\nu}^{\mathrm{NCDT}}:=\mathcal{Z}_{\sigma,\lambda,\nu,\mathrm{id}_{\Zh}}.
\]
\end{defn}
\begin{rem}\label{rem_poly}
Note that 
\[
\mathcal{Z}_{\sigma,\lambda,\nu}^{\mathrm{NCDT}}\cdot w_{\plt}\Bigl(D^{\mathrm{max}}_{\sigma,\lambda,\nu.\mathrm{id}}\Bigr)^{-1}
\]
is a formal power series in $q_+,q_-$ and $q_1,\ldots,q_{L-1}$.
\end{rem}

\begin{NB}
\section{Moduli spaces}\label{sec_moduli}
\subsection{Generating function via moduli spaces}
\subsubsection{}
Let $M^{\mathrm{max}}$ be the $A$-module given in \S \ref{} and $\{b[D^{\mathrm{max}};x,y,z]\}$ be the basis given in \S \ref{}.
\begin{lem}
There exists a finite subset $S\subset \mathcal{F}$ such that the $A$-module homomorphism given by 
\begin{equation}\label{eq_presentation}
\begin{array}{ccc}
A^{|S|} & \to & M^{\mathrm{max}},\\
\mathrm{id}^{(x,y)} & \mapsto & b[D^{\mathrm{max}};x,y,0],
\end{array}
\end{equation}
is surjective and the kernel is finitely generated, where $\mathrm{id}^{(x,y)}$ is the identity element in the component corresponding to $(x,y)\in S$.
\end{lem}

\subsubsection{}
We consider a new quiver $Q_{\sigma,\lambda,\nu,\theta}$ given by adding a new vertex $\infty$ and new edges $\{\iota_{(n,m)}\mid (n,m)\in S\}$ to the quiver $Q_{\sigma,\theta}$, where $\iota_{(n,m)}$ is an edge from $\infty$ to $i$ for $(n,m)\in \mathcal{F}_i\cap S$.
Take a generator 
\[
\Biggl\{\sum_{(x,y)\in S} a_{(x,y)}^r\cdot b[D^{\mathrm{max}};x,y,0]\Biggr\}_r\quad (a_{(x,y)}^r\in A)
\]
of the kernel of the map \eqref{eq_presentation}.
Then we consider new relations given by the original ones plus 
\[
\sum_{(x,y)\in S} a_{(x,y)}^r\circ \iota_{(x,y)}
\]
for each $r$. 
Let $A_{\sigma,\lambda,\nu,\theta}$ denote the new quiver with relations. 

The lemma below is clear from the construction:
\begin{lem}
The following are equivalent:
\begin{itemize}
\item giving a finite dimensional quotient $A$-module $M$ of $M^{\mathrm{max}}$, and 
\item giving a finite dimensional $A_{\sigma,\lambda,\nu,\theta}$-module $\tilde{M}$ with $\dim \tilde{M}_\infty=1$.
\end{itemize}
\end{lem}
\begin{prop}
For $\mathbf{v}\in \Z^I$, the moduli space $\mathcal{M}_{\sigma,\lambda,\nu,\theta}(\mathbf{v})$ exists.
\end{prop}
\begin{proof}

\end{proof}
\begin{prop}
\end{prop}

\subsubsection{}
\begin{prop}
\end{prop}
\begin{proof}

\end{proof}

\subsubsection{}
Here, we will give an explicit description of a quiver in the case $\nu=(\emptyset,\emptyset)$.

For integers $r\neq p$, we define an element $h_{[r,p]}\in A$ by 
\[
h_{[r,p]}:=
\begin{cases}
h^+_{p-1/2}\circ\cdots\circ h^+_{r+1/2},& (r<p),\\
h^-_{p+1/2}\circ\cdots\circ h^-_{r-1/2}, ,& (r>p).
\end{cases}
\]
We take the sequence of integers
\[
r_0<p_{1/2}<r_1<\cdots<p_{K-1/2}<r_K
\]
such that 
\begin{align*}
P_{\sigma,\lambda,\theta}&=
\{n\in \Z\mid G_{\sigma,\lambda,\theta}(n)>G_{\sigma,\lambda,\theta}(n\pm 1)\}=\{p_{1/2},\ldots,p_{K-1/2}\},\\
R_{\sigma,\lambda,\theta}&=
\{n\in \Z\mid G_{\sigma,\lambda,\theta}(n)<G_{\sigma,\lambda,\theta}(n\pm 1)\}=\{r_{0},\ldots,r_{K}\}.
\end{align*}
Then we have the following presentation of the $A$-module $M^{\mathrm{max}}$:
\[
\bigoplus_{k=1/2}^{K-1/2}A_{(p_k)}\to \bigoplus_{k=0}^{K} A_{(r_k)}\to M^{\mathrm{max}}\to 0.
\]
Here $A_{(p_k)}$ (resp. $A_{(r_k)}$) is isomorphic to the projective module $A_i$ if $p_k\in \pi^{-1}(i)$ (resp. $r_k\in \pi^{-1}(i)$). 
The first map is given by 
\[
e_{(p_k)}\mapsto (0,0,\ldots,0,h_{[r_{k-1/2},p_k]},-h_{[r_{k+1/2},p_k]},0,\ldots)
\]
and the second map is given by
\[
e_{(r_k)}\mapsto b[D^{\mathrm{max}};r_k,F_{\sigma,\lambda,\theta},0],
\]
where $e_{(p_k)}$ or $e_{(r_k)}$ is the idempotent in $A_{(p_k)}\subset A$ or $A_{(r_k)}\subset A$.
This gives an explicit quiver with relations.
\begin{ex}
In the case as in Example \ref{}, we have
\[
P_{\sigma,\lambda,\mathrm{id}}=\{-2\},\quad
R_{\sigma,\lambda,\mathrm{id}}=\{-3,1\}.
\]
The new quiver is given as in Figure \ref{fig28} and the new relation
\[
h^-_{-3/2}\circ h^-_{-1/2}\circ h^-_{1/2}\circ \iota_1=
h^+_{-5/2}\circ \iota_{-3}.
\]
\begin{figure}[htbp]
  \centering
  \input{fig28.tpc}
  \caption{}
\label{fig28}
\end{figure}
\end{ex}

\subsection{Potential and symmetric obstruction theory}
Assume $\nu=(\emptyset,\emptyset)$.
We consider a newer quiver $\tilde{Q}_{\sigma,\lambda,\theta}$ given by adding new edges $\{\tau_{r}\mid r\in R_{\sigma,\lambda,\theta}\}$ to the quiver $Q_{\sigma,\lambda,\vec{\emptyset},\theta}$, where $\tau_{r}$ is an edge from $i$ to $\infty$ for $r\in\pi^{-1}(i)$. 
We define a potential for this quiver by
\[
w_{\sigma,\theta}+\sum_{k=1/2}^{K-1/2}\Bigl(\tau_{p_k}\circ h_{[r_{k-1/2},p_k]}\circ\iota_{k-1/2}-\tau_k\circ h_{[r_{k+1/2},p_k]}\circ\iota_{r_{k+1/2}}\Bigr),
\]
where $w_{\sigma,\theta}$ is the original potential.

We can show the following claims just as in \cite[\S]{nagao-nakajima}. 
In fact, the claims in \cite[\S]{nagao-nakajima} are the 
\begin{ex}
The quiver is given as in Figure \ref{fig30} and the new potential is given by 
\[
w_{\sigma,\theta}+\tau_{-2}\circ h^+_{-5/2}\circ\iota_{-3}-\tau_{-2}\circ h^-_{3/2}\circ h^-_{-1/2}\circ h^-_{1/2}\circ\iota_{1}.
\]
\begin{figure}[htbp]
  \centering
  \input{fig30.tpc}
  \caption{}
\label{fig30}
\end{figure}
\end{ex}
\end{NB}

\section{Dimer shuffling and ``wall-crossing" formula}\label{sec_shuffling_wc}

\subsection{Dimer shuffling at a hexagon}\label{subsec_shuffling_6}
In this and next subsections, we study the relation between dimer configurations of type $(\plmt)$ and of type $(\sigma,\lambda,\nu,\mu_i(\theta))$ for $i\in I_H(\sigma,\theta)$.

\subsubsection{}
For $(n,m)\in\mathcal{F}$ and $M\in \Z_{>0}\sqcup \{\infty\}$. we put
\[
f(n,m;\pm,M):=\bigcup_{m'=0}^{M-1}f(n,m\pm m')
\]
We define $\partial f(n,m;\pm,M)$ and $\partial^\pm f(n,m;\pm,M)$ in the same way as in \S \ref{subsub_217} and \S \ref{subsub_322}.

\subsubsection{}
For a dimer configuration $D$ and $n\in B^{i,\pm}_{\plt}$, let $m(D,n)$ denote the unique integer such that
\[
\partial f(n,m(D,n);\sigma(i),\infty)\cap D=\partial^\pm f(n,m(D,n);\sigma(i),\infty).
\]

\subsubsection{}
For a dimer configuration $D$ and $i\in I$, we consider the following conditions:
\begin{align}
&\text{$\partial f\cap D\neq \partial^- f$ for any $f\in \mathcal{F}_i$},\label{eq_cond}\\
&\text{$\partial f\cap D\neq \partial^+ f$ if $f\in \mathcal{F}_i\backslash\{f(n,m(D,n))\mid n\in B^{i,\pm}_{\plt}\}$},\label{eq_cond2}\\
&\text{$\partial f(n,m(D,n)-2\sigma(i))\cap D\neq \partial^-f(n,m(D,n)-2\sigma(i))$ for $n\in B^{i,\pm}_{\plt}$.}\label{eq_cond3}
\end{align}

\subsubsection{}
For a dimer configuration $D^\circ$ of type $(\plmt)$
satisfying the condition \eqref{eq_cond},
we set
\[
E_i(D^\circ):=\bigl\{(n,m)\in\mathcal{F}_i\mid \partial f(n,m)\cap D^\circ=\partial^+ f(n,m)\bigr\},
\]
and define the map $M^i_{D^\circ}\colon E_i(D^\circ)\to \Z_{>0}\sqcup \{\infty\}$ by
\[
M^i_{D^\circ}(n,m):=\max\{M\mid \partial f(n,m;\sigma(i),M)\cap D^\circ=\partial^+ f(n,m;\sigma(i),M)\}.
\]
Note that
\[
\bigl(M^i_{D^\circ}\bigr)^{-1}(\infty)=\bigl\{(n,m_n)\mid n\in B^{i,+}_{\plt}\bigr\}.
\]
We put $E_i^{\mathrm{fin}}(D^\circ):=E_i(D^\circ)\backslash \bigl(M^i_{D^\circ}\bigr)^{-1}(\infty)$.

\subsubsection{}
\begin{defn}\label{defn_mutation}
For a dimer configuration $D^\circ$ of type $(\plmt)$ 
satisfying the condition \eqref{eq_cond}, let $\mu_i(D)$ be the a dimer configuration of type $(\sigma,\lambda,\nu,\mu_i(\theta))$ given by
\begin{align*}
\mu_i(D^\circ):=&\,\Biggl(D^\circ\Big\backslash\Biggl(\,\bigcup_{(n,m)\in E_i(D^\circ)}\partial^+f\bigl(n,m;\sigma(i),M^i_{D^\circ}(n,m)\bigr)\cup\bigcup_{n\in B^{i,-}_{\plt}}\partial^-f(n,m;\sigma(i),\infty)\Biggr)\Biggr)\\
&\quad\quad\quad\quad\sqcup \Biggl(\,\bigcup_{(n,m)\in E_i(D^\circ)}\partial^-f\bigl(n,m;\sigma(i),M^i_{D^\circ}(n,m)\bigr)\cup\bigcup_{n\in B^{i,-}_{\plt}}\partial^+f(n,m;\sigma(i),\infty)\Biggr).
\end{align*}
\end{defn}
Note that $\mu_i(D^\circ)$ satisfies the condition \eqref{eq_cond2} and \eqref{eq_cond3}.
\begin{ex}
In Figure \ref{fig25}, we show some examples of dimer shuffling at hexagons.
\begin{figure}[htbp]
  \centering
  \input{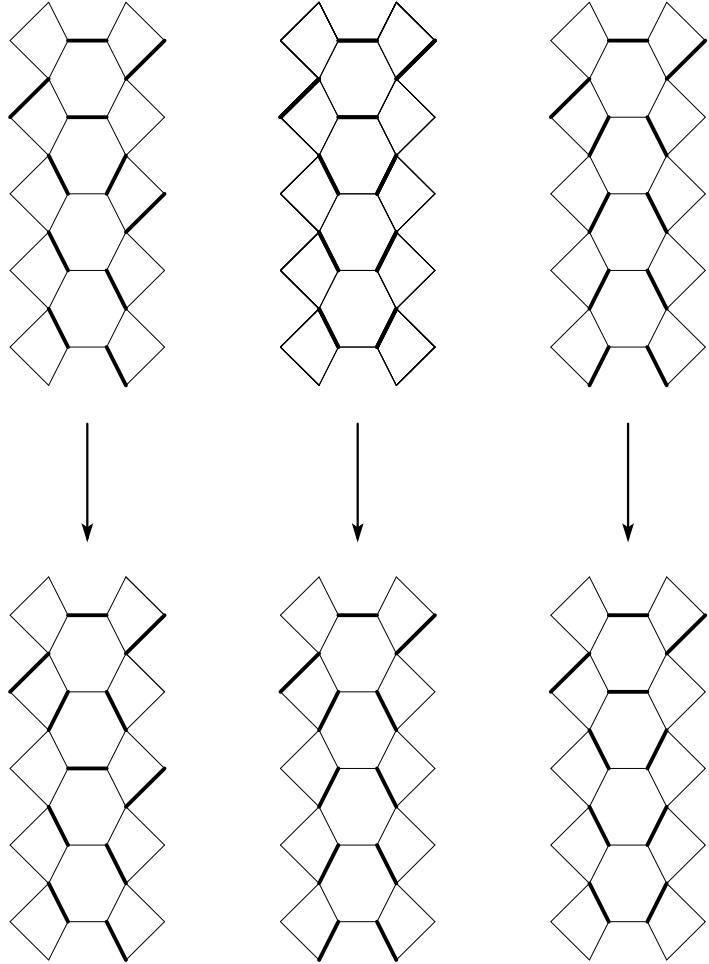}
  \caption{Examples of dimer shuffling at hexagons}
\label{fig25}
\end{figure}
\end{ex}
\begin{lem}
\[
w_{\sigma,\lambda,\mu_i(\theta)}(\mu_i(D^\circ))=w_{\plt}(D^\circ).
\]
\end{lem}
\begin{proof}
For $n\in \pi^{-1}(i)$ and $m\in\Z$ such that $n+m$ is odd, we put
\[
D^\circ(n,m):=\{\es(n+\varepsilon_1,m+\varepsilon_2)\,(\varepsilon_1,\varepsilon_2=\pm1/2)\}\cap D^\circ.
\]
Assume that
\begin{equation}\label{eq_411}
(n,m-1),(n,m+1)\notin \bigcup_{(n,m)\in E_i(D^\circ)}f\bigl(n,m;\sigma(i),M^i_{D^\circ}(n,m)\bigr),
\end{equation}
then $D^\circ(n,m)$ is one of the following:
\[
\emptyset,\,
\{\es(n\pm 1/2,m\pm 1/2)\},\,
\{\es(n\pm 1/2,m\mp 1/2)\}.
\]
In particular, we have
\[
w_{\plt}(D^\circ(n,m))=w_{\sigma,\lambda,\mu_i(\theta)}(D^\circ(n,m)).
\]
Hence we have 
\[
w_{\plt}\bigl(D^\circ\cap \mu_i(D^\circ)\bigr)=
w_{\sigma,\lambda,\mu_i(\theta)}\bigl(D^\circ\cap \mu_i(D^\circ)\bigr).
\]
The claim follows this and 
\[
w_{\plt}\Bigl(\partial^\pm f(n,m,M)\Bigr)=
w_{\sigma,\lambda,\mu_i(\theta)}\Bigl(\partial^\mp f(n,m,M)\Bigr)
\]
for $n\in \pi^{-1}(i)$.
\end{proof}

\subsection{``Wall-crossing" formula at a hexagon}\label{subsec_wc_6}
\subsubsection{}
\begin{lem}\label{lem_43}
\begin{align*}
\mathcal{Z}_{\plmt}
&=
\sum_{D^\circ}w_{\plt}(D^\circ)\cdot 
\prod_{n\in B^{i,+}_{\plt}}\bigl(1+w_{\plt}(n)\bigr)^{-1}\\
& \quad \quad \quad \quad 
\Biggl(
\prod_{(n,m)\in E^{\mathrm{fin}}_i(D^\circ)}\Bigl(1+w_{\plt}(n)^{M^i_{D^\circ}(n,m)+1}\Bigr)\big/\bigl(1+w_{\plt}(n)\bigr)\Biggr),
\end{align*}
where the sum is taken over all dimer configurations $D^\circ$ of type $(\plmt)$ satisfying the condition \eqref{eq_cond}. 
\end{lem}
\begin{proof}
For a map $s\colon E_i(D^\circ)\to \Z_{\geq 0}$ such that $s(n,m)\leq M^i_{D^\circ}(n,m)$, we define the dimer configuration $D^\circ_s$ by 
\[
D^\circ_s:=
\,\Biggl(D^\circ\Big\backslash\,\bigcup_{(n,m)\in E_i(D^\circ)}\partial^+f(n,m;\sigma(i),s(n,m))\Biggr)\sqcup\bigcup_{(n,m)\in E_i(D^\circ)}\partial^-f(n,m;\sigma(i),s(n,m)).
\]
Then we have
\[
w_{\plt}(D^\circ_s)=
w_{\plt}(D^\circ)\cdot 
\prod_{(n,m)\in E_i(D^\circ)}w_{\plt}(n)^{s(n,m)}.
\]

Note that any dimer configuration $D$ is realized as $D^\circ(s)$ by some $D^\circ$ and $s$ uniquely.
Hence we have 
\begin{align*}
\mathcal{Z}_{\plmt}
&=
\sum_{D^\circ}w_{\plt}(D^\circ)\cdot \biggl(\sum_s\prod_{(n,m)\in E_i(D^\circ)}w_{\plt}(n)^{s(n,m)}\biggr)\\
&=
\sum_{D^\circ}w_{\plmt}(D^\circ)\cdot 
\prod_{n\in B^{i,+}_{\plt}}\bigl(1-w_{\plt}(n)\bigr)^{-1}\\
& \quad \quad \quad \quad 
\Biggl(\,
\prod_{(n,m)\in E_i^{\mathrm{fin}}(D^\circ)}\Bigl(1+w_{\plt}(n)^{M^i_{D^\circ}(n,m)+1}\bigr)\big/\bigl(1+w_{\plt}(n)\Bigr)\Biggr).
\end{align*}
\end{proof}

\subsubsection{}
\begin{thm}\label{thm_hex}
\[
\mathcal{Z}_{\sigma,\lambda,\nu,\mu_i(\theta)}
=\mathcal{Z}_{\sigma,\lambda,\nu,\theta}
\cdot
\prod_{n\in B^{i,+}_{\sigma,\lambda,\theta}}\bigl(1-w_{\sigma,\lambda,\theta}(n)\bigr)
\cdot
\prod_{n\in B^{i,-}_{\sigma,\lambda,\theta}}\bigl(1-w_{\sigma,\lambda,\theta}(n)\bigr)^{-1}.
\]
\end{thm}
\begin{proof}
As Lemma \ref{lem_43}, we get
\begin{align*}
\mathcal{Z}_{\sigma,\lambda,\nu,\mu_i(\theta)}
&=
\sum_{D^\bullet}w_{\sigma,\lambda,\mu_i(\theta)}(D^\bullet)\cdot 
\prod_{n\in B^{i,+}_{\sigma,\lambda,\mu_i(\theta)}}\Bigl(1-w_{\sigma,\lambda,\mu_i(\theta)}(n)^{-1}\Bigr)^{-1}\\
& \quad 
\Biggl(\,
\prod_{(n,m)\in \check{E}_i(D^\bullet)}\Bigl(1+w_{\sigma,\lambda,\mu_i(\theta)}(n)^{-\check{M}^i_{D^\bullet}(n,m)-1}\Bigr)\Big/\Bigl(1+w_{\sigma,\lambda,\mu_i(\theta)}(n)^{-1}\Bigr)\Biggr),
\end{align*}
where the sum is taken over all dimer configurations $D^\bullet$ of type $(\sigma,\lambda,\nu,\mu_i(\theta))$ satisfying the condition \eqref{eq_cond2} and \eqref{eq_cond3} and
\begin{align*}
\check{E}_i(D^\bullet)&:=\bigl\{(n,m)\in\mathcal{F}_i\,\big|\,\partial f(n,m)\cap D^\bullet=\partial^-f(n,m)\bigr\},\\
\check{M}^i_{D^\bullet}(n,m)&:=\max\bigl\{M\,\big|\,\partial f(n,m;\sigma(i),M)\cap D^\bullet=\partial^- f(n,m;-\sigma(i),M)\bigr\}.
\end{align*}
Note that $\mu_i$ gives a one-to-one correspondence of dimer configurations of type $(\plmt)$ satisfying the condition \eqref{eq_cond} and ones of type $(\sigma,\lambda,\nu,\mu_i(\theta))$ satisfying the condition \eqref{eq_cond2} and \eqref{eq_cond3}.
Hence the claim follows from
\begin{itemize}
\item $B^{i,\pm}_{\sigma,\lambda,\mu_i(\theta)}=B^{i,\mp}_{\sigma,\lambda,\theta}$,
\item $w_{\sigma,\lambda,\mu_i(\theta)}(n)=w_{\sigma,\lambda,\theta}(n)^{-1}$ for $n\in \pi^{-1}(i)$,
\item $(n,m)\mapsto (n,m+\sigma(i)\cdot(M^i_{D^\circ}(n,m)-1))$ gives a bijection between $E_i^{\mathrm{fin}}(D^\circ)$ and $\check{E}_i(\mu_i(D^\circ))$ which respects $M^i_{D^\circ}$ and $\check{M}^i_{\mu_i(D^\circ)}$, 
\end{itemize}
and Lemma \ref{lem_43}.
\end{proof}

\subsection{Dimer shuffling at a quadrilateral}\label{subsec_shuffling_4}
In this and next subsections, we study the relation between dimer configurations of type $(\plmt)$ and of type $(\sigma,\lambda,\nu,\mu_i(\theta))$ for $i\in I_S(\sigma,\theta)$.

\subsubsection{}
For a dimer configuration $D^\circ$ of type $(\plmt)$ satisfying the condition \eqref{eq_cond} and $n\in\pi^{-1}(i)$, we define
\begin{align*}
E^1_n(D^\circ)&:=\{(n,m)\in\mathcal{F}\mid \partial f(n,m)\cap D^\circ=\partial^+ f(n,m)\},\\
E^2_n(D^\circ)&:=\{(n,m)\in\mathcal{F}\mid \partial f(n,m)\cap D^\circ=\emptyset\,\}.
\end{align*}
\begin{lem}\label{lem_45}
\[
|E^1_n(D^\circ)|-|E^2_n(D^\circ)|=
\begin{cases}
\mp 1, & n\in B^{i,\pm}_{\plt},\\
0, & \text{otherwise}.
\end{cases}
\]
\textup{(}See \eqref{eq_Bpm2} for the notation.\textup{)}
\end{lem}
\begin{proof}
For $n,m\in \Z$ such that $n+m$ is odd, we define $\varepsilon_{D^\circ}(n,m)$  by
\[
\varepsilon_{D^\circ}(n,m):=\begin{cases}
+, & \es(n+1/2,m+1/2),\es(n-1/2,m+1/2)\notin D,\\
-, & \es(n+1/2,m-1/2),\es(n-1/2,m-1/2)\notin D.\\
\end{cases}
\]
Then for $(n,m)\in\mathcal{F}$, we have
\begin{align*}
(n,m)\in E^1_n(D^\circ) &\iff \varepsilon_{D^\circ}(n,m\pm 1)=\pm,\\
(n,m)\in E^2_n(D^\circ) &\iff \varepsilon_{D^\circ}(n,m\pm 1)=\mp,
\end{align*}
and 
$\varepsilon_{D^\circ}(n,m)=\mp\tlam(n\pm 1/2)$ if $\tsig(n\pm 1/2)\cdot m \gg 0$.
Thus the claim follows.
\end{proof}

\subsubsection{}\label{subsub_432}
For a dimer configuration $D^\circ$ of type $(\plmt)$ satisfying the condition \eqref{eq_cond}, we define a new dimer configuration $\mu_i(D^\circ)$ of type $(\sigma,\lambda,\nu,\mu_i(\theta))$ as follows:
\begin{itemize}
\item if $\pi(h)\neq i\pm 1/2$ then we have 
\[
\es(h,k)\in D^\circ\iff\es(h,k)\in \mu_i(D^\circ),
\]
\item if $n\in I_H(\sigma,\theta)$ and $\pi(n)\neq i\pm 1$ then we have 
\[
\eh(n,m)\in D^\circ\iff\eh(n,m)\in \mu_i(D^\circ),
\]
\item for $(n,m)\in\mathcal{F}_i$ we have
\begin{align*}
D^\circ(f(n,m))=\emptyset&\iff \mu_i(D^\circ)(f(n,m))=\partial^-_{\sigma,\mu_i(\theta)}(f(n,m)),\\
D^\circ(f(n,m))=\partial^+_{\sigma,\theta}(f(n,m))&\iff \mu_i(D^\circ)(f(n,m))=\emptyset,
\end{align*}
\end{itemize}
\begin{rem}
Here we use such notations as $\partial^\pm_{\sigma,\theta}(f(n,m))$, in order to emphasize that the notions like $\partial^\pm(f(n,m))$ given in \S \ref{subsub_217} depend on $\sigma$ and $\theta$.
\end{rem}
\begin{itemize}
\item if $D^\circ(f(n,m))\neq \emptyset,\partial^+_{\sigma,\theta}(f(n,m))$ for $(n,m)\in\mathcal{F}_i$, then we have
\[
\es(n+\varepsilon_1,m+\varepsilon_2)\in D^\circ\iff
\es(n-\varepsilon_1,m-\varepsilon_2)\in \mu_i(D^\circ)\quad(\varepsilon_1,\varepsilon_2=\pm1/2),
\]
\item if $\sigma(i\pm 3/2)\neq\sigma(i\pm 1/2)$ then we have
\[
\es(n\pm 1/2,m-1),\es(n\pm 1/2,m+1)\notin D^\circ\iff \eh(n\pm 1,m)\in \mu_i(D^\circ).
\]
\end{itemize}
Note that $\mu_i(D^\circ)$ satisfies the following condition:.
\begin{equation}\label{eq_cond4}
\text{$D(f)\neq \partial^+ f$ for any $f\in \mathcal{F}_i$}.
\end{equation}
\begin{ex}
In Figure \ref{fig26}, we show some examples of dimer shuffling at squares.
\begin{figure}[htbp]
  \centering
  \input{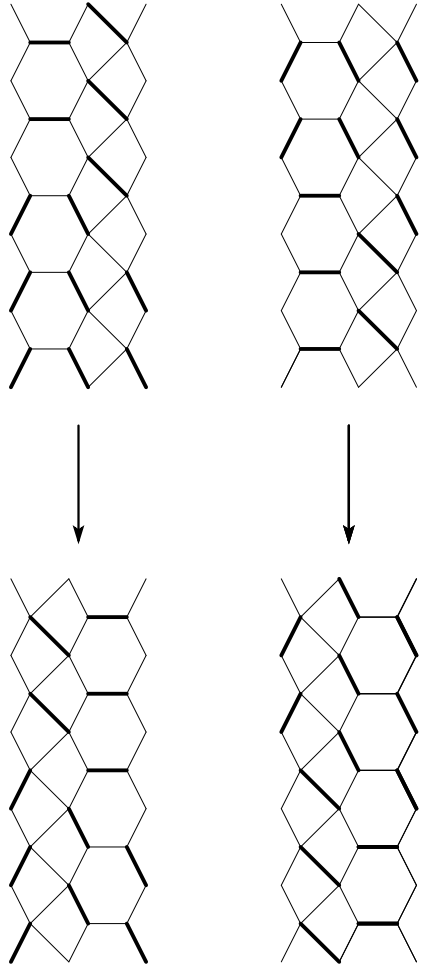}
  \caption{Examples of dimer shuffling at squares}
\label{fig26}
\end{figure}
\end{ex}

\subsubsection{}
\begin{lem}
\[
w_{\sigma,\lambda,\mu_i(\theta)}(\mu_i(D^\circ))=
w_{\plt}(D^\circ).
\]
\end{lem}
\begin{proof}
Note that we have 
\[
w_{\plt}(\partial^+_{\sigma,\theta}f)=
w_{\sigma,\lambda,\mu_i(\theta)}(\partial^-_{\sigma,\mu_i(\theta)}f)
\]
for $f\in\mathcal{F}_i$ and 
\[
w_{\plt}(\partial^+_{\sigma,\theta}f)=
\begin{cases}
1, & n\in B^{i,+}_{\sigma,\lambda,\theta},\\
w_{\plt}(n)^{-1}, & n\in B^{i,-}_{\sigma,\lambda,\theta}.
\end{cases}
\]
Thus, the claim follows from Lemma \ref{lem_45} and \eqref{eq_def_wt}.
\end{proof}

\subsection{``Wall-crossing" formula at a quadrilateral}\label{subsec_wc_4}
\subsubsection{}
\begin{lem}\label{lem_48}
\[
\mathcal{Z}_{\plmt}
=
\sum_{D^\circ}w_{\plt}(D^\circ)\cdot 
\prod_{n\in \pi^{-1}(i)}\bigl(1+w_{\plt}(n)\bigr)^{|E_n^1(D^\circ)|}.
\]
\end{lem}
\begin{proof}
We set 
\[
E^1_i(D^\circ):=\bigcup_{n\in \pi^{-1}(i)}E^1_n(D^\circ),\quad
E^2_i(D^\circ):=\bigcup_{n\in \pi^{-1}(i)}E^2_n(D^\circ).
\]
Given a subset $S\subset E^1_i(D^\circ)$, we get the dimer configuration $D^\circ_S$ of type $(\plmt)$ such that
\[
D^\circ_S:=\Bigl(D\backslash \bigcup\partial^+f\Bigr)\cup\bigcup\partial^+f
\]
and we have
\[
w_{\plt}(D^\circ_S)=w_{\plt}(D^\circ)\cdot \prod_{(n,m)\in S}w_{\plt}(n).
\]
Note that any dimer configuration $D$ is realized as $D^\circ_S$ by some $D^\circ$ and $S$ uniquely.
Hence we have 
\begin{align*}
\mathcal{Z}_{\plmt}
&=
\sum_{D^\circ}w_{\plt}(D^\circ)\cdot \biggl(\sum_S\prod_{(n,m)\in S}w_{\plt}(n)\biggr)\\
&=
\sum_{D^\circ}w_{\plt}(D^\circ)\cdot 
\prod_{(n,m)\in E^1_i(D^\circ)}\bigl(1+w_{\plt}(n)\bigr)\\
&=
\sum_{D^\circ}w_{\plt}(D^\circ)\cdot 
\prod_{n\in \pi^{-1}(i)}\bigl(1+w_{\plt}(n)\bigr)^{|E_n^1(D^\circ)|}.
\end{align*}
\end{proof}

\subsubsection{}
\begin{thm}\label{thm_quad}
\[
\mathcal{Z}_{\sigma,\lambda,\nu,\mu_i(\theta)}=
\mathcal{Z}_{\sigma,\lambda,\nu,\theta}
\cdot
\prod_{n\in B^{i,+}_{\sigma,\lambda,\theta}}\bigl(1+w_{\sigma,\lambda,\theta}(n)\bigr)^{-1}
\cdot
\prod_{n\in B^{i,-}_{\sigma,\lambda,\theta}}\bigl(1+w_{\sigma,\lambda,\theta}(n)\bigr).
\]
\end{thm}
\begin{proof}
Let $D^\bullet$ be a dimer configuration of type $(\sigma,\lambda,\nu,\mu_i(\theta))$ satisfying the condition \eqref{eq_cond4}. 
We put
\[
\tilde{E}^1_n(D^\bullet):=\bigr\{(n,m)\in\mathcal{F}\,\big|\, \partial_{\sigma,\mu_i(\theta)} f(n,m)\cap D^\bullet=\partial^-_{\sigma,\mu_i(\theta)} f(n,m)\big\},
\]
then, as Lemma \ref{lem_48}, we get
\[
\mathcal{Z}_{\sigma,\lambda,\nu,\mu_i(\theta)}
=
\sum_{D^\bullet}w_{\sigma,\lambda,\mu_i(\theta)}(D^\bullet)\cdot 
\prod_{n\in \pi^{-1}(i)}\bigl(1+w_{\sigma,\lambda,\mu_i(\theta)}(n)^{-1}\bigr)^{|\tilde{E}_n^1(D^\bullet)|},
\]
where the sum is taken over all dimer configurations $D^\bullet$ of type $(\sigma,\lambda,\nu,\mu_i(\theta))$ satisfying the condition \eqref{eq_cond4}. 
Note that $\mu_i$ gives a one-to-one correspondence of dimer configurations of type $(\plmt)$ satisfying the condition \eqref{eq_cond} and ones of type $(\sigma,\lambda,\nu,\mu_i(\theta))$ satisfying the condition \eqref{eq_cond4}.
Hence the claim follows from 
\begin{itemize}
\item $\tilde{E}_n^1(\mu_i(D^\circ))=E_n^2(D^\circ)$ for $n\in \pi^{-1}(i)$, 
\item $w_{\sigma,\lambda,\mu_i(\theta)}(n)=w_{\sigma,\lambda,\theta}(n)^{-1}$ for $n\in \pi^{-1}(i)$
\end{itemize}
and Lemma \ref{lem_48}.
\end{proof}

\subsection{Conclusion}\label{subsec_46}
\subsubsection{}\label{thm_410}
For $\sigma$ and $\alpha\in \Lambda^{\mathrm{re},+}$, we put 
\begin{equation}\label{eq_sigma_alpha}
\sigma(\alpha):=\sigma(j^-(\alpha))\cdot\sigma(j^+(\alpha)).
\end{equation}
Combining Theorem \ref{thm_hex} and \ref{thm_quad}, we get the following:
\begin{thm}
\begin{align*}
\mathcal{Z}_{\sigma,\lambda,\nu,\theta}&=
\mathcal{Z}_{\sigma,\lambda,\nu}^{\mathrm{NCDT}}
\cdot
\prod_{\alpha\in\Lambda_\theta^{\mathrm{re},+}}\Biggr(\prod_{(h,h')\in B^{\alpha,+}_{\sigma,\lambda}}(1-\sigma(\alpha)\cdot w_\lambda(h')/w_\lambda(h))^{\sigma(\alpha)}\\
&\hspace{50mm}\prod_{(h,h')\in B^{\alpha,-}_{\sigma,\lambda}}(1-\sigma(\alpha)\cdot w_\lambda(h')/w_\lambda(h))^{-\sigma(\alpha)}\Biggl).
\end{align*}
\textup{(}See \eqref{eq_lam_theta} and \eqref{eq_Bpm1} for the notations.\textup{)}
\end{thm}
Since the second term in the right-hand side does not depend on $\nu$, we have the following:
\begin{cor}\label{cor_413}
\[
\mathcal{Z}_{\sigma,\lambda,\nu,\theta}\big/\mathcal{Z}_{\sigma,\lambda,\vec{\emptyset},\theta}=
\mathcal{Z}_{\sigma,\lambda,\nu}^{\mathrm{NCDT}}\big/\mathcal{Z}_{\sigma,\lambda,\vec{\emptyset}}^{\mathrm{NCDT}}.
\]
\end{cor}

\subsubsection{}
By Lemma \ref{lem_14} and Theorem \ref{thm_410}, we have the following:
\begin{thm}\label{thm_412}
\begin{align*}
&\mathcal{Z}_{\sigma,\lambda,\nu,\theta}\Big|_{q_+=q_-=(q_0)^{1/2}}=\\
&\quad \mathcal{Z}_{\sigma,\lambda,\nu}^{\mathrm{NCDT}}\Big|_{q_+=q_-=(q_0)^{1/2}}
\cdot
\prod_{
\alpha\in \Lambda_\theta^{\mathrm{re},+}}
\bigl(1-\sigma(\alpha)\cdot q^{\alpha}\bigr)^{\sigma(\alpha)\cdot \bigl(\alpha^0+c_\lambda[j_-(\alpha)]-c_\lambda[j_+(\alpha)]\bigr)}.
\end{align*}
\textup{(}See \eqref{eq_lam_theta} for the notation.\textup{)}
\end{thm}
Since the second term in the right-hand side depend only on $c_\lambda[j]$'s but not on $\lambda$ and $\nu$, we have the following:
\begin{cor}\label{cor_415}
If $c_\lambda[j]=0$ for any $j$, then we have
\[
\mathcal{Z}_{\sigma,\lambda,\nu,\theta}\big/\mathcal{Z}_{\sigma,\vec{\emptyset},\vec{\emptyset},\theta}\Big|_{q_+=q_-}
=
\mathcal{Z}_{\sigma,\lambda,\nu}^{\mathrm{NCDT}}\big/\mathcal{Z}_{\sigma,\vec{\emptyset},\vec{\emptyset}}^{\mathrm{NCDT}}\Big|_{q_+=q_-}
.
\]
\end{cor}

\section{Refined topological vertex via dimer model}\label{sec_rtv}
\subsection{Refined topological vertex for $\C^3$}\label{subsec_rtv}

\subsubsection{}
Note that a Young diagram can be regarded as a subset of $(\Z_{\geq 0})^2$.
For a Young diagram $\lambda$, let $\Lambda^x(\lambda)$ (resp. $\Lambda^y(\lambda)$ or $\Lambda^z(\lambda)$) be the subset of $(\Z_{\geq 0})^3$ consisting of the elements $(x,y,z)\in(\Z_{\geq 0})^3$ such that 
$(y,z)\in\lambda$ (resp. $(z,x)\in\lambda$ or $(x,y)\in\lambda$).

\subsubsection{}
Given a triple $(\lambda_x, \lambda_y, \lambda_z)$ of Young diagrams, let $\Lambda^{\mathrm{min}}=\Lambda^{\mathrm{min}}_{\lambda_x, \lambda_y, \lambda_z}$ be the following subset of $(\Z_{\geq 0})^3$:
\[
\Lambda^{\mathrm{min}}:=\Lambda^x(\lambda_x)\cup\Lambda^y(\lambda_y)\cup\Lambda^z(\lambda_z)\subset(\Z_{\geq 0})^3.
\]
\subsubsection{}
A subset $\Lambda$ of $(\Z_{\geq 0})^3$ is said to be a {\it $3$-dimensional Young diagram of type $(\lambda_x, \lambda_y, \lambda_z)$} if the following conditions are satisfied:
\begin{itemize}
\item if $(x,y,z)\notin\Lambda$, then $(x+1,y,z), (x,y+1,z), (x,y,z+1)\notin\Lambda$,
\item $\Lambda\supset\Lambda^{\mathrm{min}}$, and
\item $|\Lambda\backslash\Lambda^{\mathrm{min}}|<\infty$.
\end{itemize}
\subsubsection{}
For a Young diagram $\lambda$, we define a monomial $w_{\lambda}(m)$ for each $m\in\Z$ by 
\begin{equation}\label{eq_refined_weight}
w_{\lambda}(m)=q_{\lambda(m-1/2)}\cdot q_{\lambda(m+1/2)}\cdot q_1\cdot\cdots\cdot q_{L-1}.
\end{equation}
For a finite subset $S$ of $(\Z_{\geq 0})^3$ we define the weight $w(S)$ by
\[
w(S):=\prod_{(x,y,z)\in S}w_{\lambda_x}(y-z).
\]
For a positive integer $N$, we set $C_N:=[0,N]^3$.
Given a $3$-dimensional Young diagram $\Lambda$ of type $(\lambda_x, \lambda_y, \lambda_z)$, we take a sufficiently large $N$ such that $\Lambda\backslash\Lambda^{\mathrm{min}}\subset C_N$ and define the weight $w(\Lambda)$ of $\Lambda$ by
\[
w(\Lambda):=w(\Lambda\cap C_N)\big/\bigl(w(\Lambda^x(\lambda_x)\cap C_N)\cdot w(\Lambda^y(\lambda_y)\cap C_N)\cdot w(\Lambda^z(\lambda_z)\cap C_N)\bigr).
\]
Note that this is well-defined.
\begin{rem}
In the definition of $w(\Lambda)$, the three axes does not play the same roles.
The $x$-axis is called the \textup{preferred axis} for the refined topological vertex. 
\end{rem}
\begin{rem}
If we replace the definition \eqref{eq_refined_weight} with 
\[
(q_{\lambda(m-1/2)})^2\cdot q_1\cdot\cdots\cdot q_{L-1},
\]
then the weight coincides with the one in \cite{RTV}. 
Our weight coincides with the one in \cite{dimofte-gukov}.
\end{rem}
We define the generating function 
\[
G_{\lambda_x, \lambda_y, \lambda_z}(\vec{q}\,):=\sum w(\Lambda)
\]
where the sum is taken over all $3$-dimensional Young diagrams of type $(\lambda_x, \lambda_y, \lambda_z)$.

\subsection{Dimer model for $L=1$}\label{subsec_51_2}
In the case $L=1$, the graph in \S \ref{subsec_graph} gives a hexagon lattice. 
As we have only two choices of $\sigma$, we put $\sigma(1/2)=+$.
We take $\mathrm{id}$ as $\theta$.
We omit $\sigma$ and $\mathrm{id}$ from the notations in this subsection. 
Note that $\lambda$ is a single $2$-dimensional Young diagram.

It is well-known that giving a dimer configuration of type $(\lambda,\nu)$ is equivalent to giving a $3$-dimensional Young diagram of type $(\lambda, \nu_+, {}^{\mathrm{t}}\nu_-)$.
Let $D(\Lambda)$ be the dimer configuration corresponding to a $3$-dimensional Young diagram $\Lambda$. 

For a Young diagram $\eta=(\eta_{(1)},\eta_{(2)},\ldots)$ and a monomial $p$, we put
\[
w(\eta;p,Q):=\prod\bigl(pQ^{i-1}\bigr)^{\eta_{(i)}}.
\]
Then we can verify the following:
\begin{equation}\label{eq_weight3}
w_{\lambda}(D(\Lambda))=w(\nu_-;q_+,Q_+)\cdot w(\nu_+;q_-,Q_-)\cdot w(\Lambda).
\end{equation}
\begin{ex}
As we show in Figure \ref{fig22}, we have
\begin{align*}
&w_\emptyset(\Lambda^{\mathrm{min}}_{\emptyset,(1),\emptyset})=w((1);q_-,Q_-)=q_-,\\ 
&w_\emptyset(\Lambda^{\mathrm{min}}_{\emptyset,(2),\emptyset})=w((2);q_-,Q_-)=q_-^2,\\
&w_\emptyset(\Lambda^{\mathrm{min}}_{\emptyset,(1,2),\emptyset})=w((2,1);q_-,Q_-)=q_-^3Q_-.
\end{align*}
\begin{figure}[htbp]
  \centering
  \input{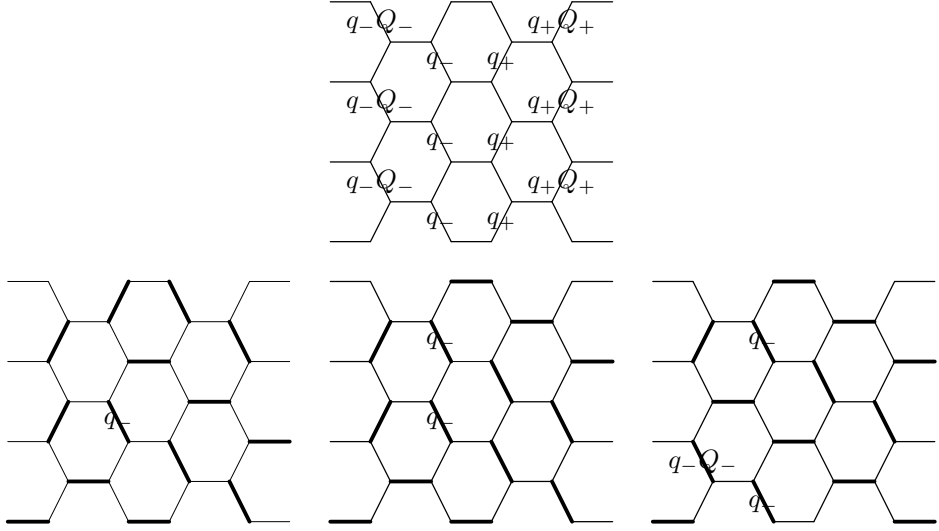}
  \caption{$D(\Lambda^{\mathrm{min}}_{\emptyset,(1),\emptyset})$, $D(\Lambda^{\mathrm{min}}_{\emptyset,(2),\emptyset})$ and $D(\Lambda^{\mathrm{min}}_{\emptyset,(1,2),\emptyset})$}
\label{fig22}
\end{figure}
\end{ex}
In particular, we have 
\[
\mathcal{Z}_{\lambda,\nu}=w(\nu_-;q_+,Q_+)\cdot w(\nu_+;q_-,Q_-)\cdot G_{\lambda, \nu_+, {}^{\mathrm{t}}\nu_-},
\]
where $\mathcal{Z}_{\lambda,\nu}$ is the generationg function given in Definition \ref{defn_gen_func}.

\subsection{Refined topological vertex for a small resolution}\label{subsub_515}
We will define generating functions $\mathcal{Z}^{\mathrm{RTV}}_{\sigma,\lambda,\nu}(\vec{q}\,)$. 
First, we consider the following data:
let $\vec{\nu}=(\nu^{(1)},\ldots,\nu^{(L-1)})$ be an ($L-1$)-tuple of Young diagrams and 
$\vec{\Lambda}=(\Lambda^{(1/2)},\ldots,\Lambda^{(L-1/2)})$ be an $L$-tuple of $3$-dimensional Young diagrams such that $\Lambda^{(j)}$ is
\begin{itemize}
\item of type $(\lambda^{(j)},\nu^{(j+1/2)},{}^{\mathrm{t}}\nu^{(j-1/2)})$ if $\sigma(j)=+$,
\item of type $(\lambda^{(j)},{}^{\mathrm{t}}\nu^{(j-1/2)},\nu^{(j+1/2)})$ if $\sigma(j)=-$,\label{item_weight}
\end{itemize}
where we put $\nu^{(0)}:=\nu_-$ and $\nu^{(L)}:=\nu_+$.
We say that the data $(\vec{\Lambda},\vec{\nu})$ is of type $(\sigma,\lambda,\nu)$.
We define the weight $w(\vec{\Lambda},\vec{\nu})$ of the data $(\vec{\Lambda},\vec{\nu})$ by 
\[
w_\sigma(\vec{\Lambda},\vec{\nu}):=w(\nu_+;q_-,Q_-)\cdot w(\nu_-;q_+,Q_+)\cdot\Biggl(\prod_{j=1/2}^{L-1/2} w(\Lambda^{(j)})\Biggr)\cdot\Biggl(\prod_{i=1}^{L-1} w^i_\sigma(\mu^{(i)})\Biggr)
\]
where $w^i_\sigma(\mu^{(i)})$ is given by
\begin{equation}\label{eq_weight2}
w^i_\sigma(\mu^{(i)}):=\prod_{(\alpha,\beta)\in\mu^i}
\begin{cases}
q_i\cdot Q^{2\alpha+1}, & \sigma(i-1/2)=\sigma(i+1/2)=+,\\
q_i\cdot Q^{2\beta+1}, & \sigma(i-1/2)=\sigma(i+1/2)=-,\\
q_i\cdot Q\cdot Q_+^{\alpha}\cdot Q_-^{\beta}, & \sigma(i-1/2)=+,\,\sigma(i+1/2)=-,\\
q_i\cdot Q\cdot Q_-^{\alpha}\cdot Q_+^{\beta}, & \sigma(i-1/2)=-,\,\sigma(i+1/2)=+.
\end{cases}
\end{equation}
We consider the following generating function:
\[
\mathcal{Z}^{\mathrm{RTV}}_{\sigma,\lambda,\mu}(\vec{q}\,):=\sum w_\sigma(\vec{\Lambda},\vec{\nu})
\]
where the sum is taken over all the data as above.
\begin{rem}
This is the generating function of the refined topological vertex associated to $Y_\sigma$, where $Y_\sigma\to X$ is the crepant resolution constructed from $\sigma$ \textup{(}see \cite[\S 1.1]{3tcy} for the construction of $Y_\sigma$\textup{)}.
In Figure \ref{fig32}, we show the polygon corresponding to $Y_\sigma$ for $\sigma$ given by 
\[
(\sigma(1/2),\ldots,\sigma(11/2))=(+,-,+,+,-,+).
\]
\begin{figure}[htbp]
  \centering
  \input{fig32.tpc}
  \caption{$Y_\sigma$}
\label{fig32}
\end{figure}
\end{rem}

\subsection{Limit behavior of the dimer model}\label{subsec_limit}
\subsubsection{}\label{subsub_521}

Let $\mathbf{i}\in I^{\Z_{>0}}$ be a minimal expression such that for any $N\in\Z_{\geq 0}$ we have $b(N)\in\Z_{>0}$ such that $\alpha_{\mathbf{i},{b}}>N\delta$ for any $b>b(N)$. 

\subsubsection{}\label{subsub_522}
\begin{lem}\label{lem_52}
Given $\sigma$, $\lambda$ and a monomial $\mathbf{q}$, 
there exists an integer  $B_1$ such that the following condition holds: 
for any $b\geq B_1$,
\begin{itemize}
\item
any dimer configuration of type $(\sigma,\lambda,\nu,\theta_{\mathbf{i},b})$ with weight $\mathbf{q}$ satisfies the condition \eqref{eq_cond}
\item
any dimer configuration of type $(\sigma,\lambda,\nu,\theta_{\mathbf{i},{b+1}})$ with weight $\mathbf{q}$ satisfies the condition \eqref{eq_cond2},
\item
$\mu_{i_b}$ gives a one-to-one correspondence between dimer configurations of type $(\sigma,\lambda,\nu,\theta_{\mathbf{i},b})$ with weight $(\sigma,\lambda,\nu,\theta_{\mathbf{i},b+1})$ with weight $\mathbf{q}$,
\end{itemize}
\end{lem}
\begin{proof}
Take $N_2$ so that 
\[
q^{N_2\delta}>\mathbf{q}\cdot w_{\plt}\Bigl(D^{\mathrm{max}}_{\sigma,\lambda,\nu.\mathrm{id}}\Bigr)^{-1}.
\]
By Remark \ref{rem_poly} and Theorem \ref{thm_412}, 
\[
\mathcal{Z}_{\sigma,\lambda,\nu,\theta}\cdot w_{\plt}\Bigl(D^{\mathrm{max}}_{\sigma,\lambda,\nu.\mathrm{id}}\Bigr)^{-1}\Big|_{q_+=q_-=(q_0)^{1/2}}
\]
is a polynomial in $q_0,\ldots,q_{L-1}$. 
Thus, there does not exist any dimer configuration with weight $\mathbf{q}-\alpha(\mathbf{i},b)$ for any $b>b(N_2)=:B_1$, where $b(N_2)$ is taken as in \S \ref{subsub_521}. 

Assume that we have a dimer configuration type $(\sigma,\lambda,\nu,\theta_{\mathbf{i},b})$ with weight $\mathbf{q}$ and $f\in \mathcal{F}$ such that $D(f)=\partial^-(f)$.
Then we get a dimer configuration $D\cup\partial^+(f) \backslash \partial^-(f)$  with weight $\mathbf{q}-\alpha(\theta,i)$, which is a contradiction. 
We can check the second claim similarly and the third claim immediately follows from the first and second ones.
\end{proof}

\subsubsection{}\label{subsub_523}
Given $\sigma$, $\lambda$, we can take an integer $N_2$ such that
\begin{itemize}
\item $\tsig(h)=\pm\tlam(h)$ for any $h\in\Zh$ such that $\pm h>N_2L$, 
\item $e^{\mathrm{s}}(h,k)\notin D^{\mathrm{max}}_{\sigma,\lambda,\theta_{\mathbf{i},B_1}}$ for any $h$ and $k$ such that $h<N_2L$ and $h\cdot\tsig(h)-k$ is even, and 
\item $e^{\mathrm{s}}(h,k)\notin D^{\mathrm{max}}_{\sigma,\lambda,\theta_{\mathbf{i},B_1}}$ for any $h$ and $k$ such that $h>N_2L$ and $h\cdot\tsig(h)-k$ is odd.
\end{itemize} 
Take a monomial $\mathbf{q}$.
Since we have only a finite number of dimer configuration of type $(\sigma,\lambda,\nu,\theta_{\mathbf{i},B_1})$ with weight $\mathbf{q}$ and each dimer configuration has only finite difference with $D^{\mathrm{max}}_{\sigma,\lambda,\nu,\theta_{\mathbf{i},B_1}}$, we can take an integer $N_4$ such that
\begin{itemize}
\item $\tsig(h)=\pm\tlam(h)$ for any $h\in\Zh$ such that $\pm h>LN_4$, 
\item $e^{\mathrm{s}}(h,k)\notin D$ for any $h$ and $k$ such that $h<LN_4$ and $h\cdot\tsig(h)-k$ is even, and 
\item $e^{\mathrm{s}}(h,k)\notin D$ for any $h$ and $k$ such that $h>LN_4$ and $h\cdot\tsig(h)-k$ is odd.
\end{itemize}

\subsubsection{}\label{subsub_524}
\begin{lem}\label{lem_53}
Let $D$ be a dimer configuration of type $(\plmt)$ satisfying the condition \eqref{eq_cond}.
Take $h\in \pi^{-1}(i+1/2)$ such that $\tsig(h)=\tlam(h)$ and assume that $\es(h,k)\notin D$ for any $k\in \Zh$ such that $h\tsig(h)-k$ is odd.
Then $\es(h-1,k-\tsig(h))\notin \mu_i(D)$.

Similarly, take $h\in \pi^{-1}(i+1/2)$ such that $\tsig(h)=-\tlam(h)$ and assume that $\es(h,k)\notin D$ for any $k\in \Zh$ such that $h\tsig(h)-k$ is even.
Then $\es(h+1,k+\tsig(h))\notin \mu_i(D)$.
\end{lem}
\begin{proof}
In the case $i\in I_S$, for any $h,k\in\Zh$ such that $\tsig(h)=\tlam(h)$ and $h\tsig(h)-k$ is odd, we can verify 
\[
\es(h,k)\notin D\Longrightarrow \es(h-1,k-\tsig(h))\notin \mu_i(D)
\]
from the definition of $\mu_i(D)$ in \S \ref{subsub_432}.

In the case $i\in I_S$, assume we have $k\in\Zh$ such that $h\tsig(h)-k$ is odd and $\es(h-1,k-\tsig(h))\in \mu_i(D)$.
From Definition \ref{defn_mutation}, we have $\es(h-1,k-\tsig(h))\in D$.
Since $\es(h,k-2\tsig(h))\notin D$, we have $\es(h,k-\tsig(h))\in D$.
Then, since $\tsig(h)=\tlam(h)$, there exists $m$ such that $\sigma(i)(m-k)>0$ and $\partial f(h-1/2,m)\cap D=\partial^- f(h-1/2,m)$, which is a contradiction.
\end{proof}

\subsubsection{}\label{subsub_525}
Given $\sigma$, $\lambda$ and a monomial $\mathbf{q}$, take $B_1$ and $N_4$ as in Lemma \ref{lem_52} and \S \ref{subsub_523}. By the definition of $N_4$ and Lemma \ref{lem_53}, we have the following lemma:

\begin{lem}\label{lem_limit}
For any $b\geq B_1$ and for any dimer configuration of type $(\sigma,\lambda,\nu,\theta_{\mathbf{i},b}$ with weight $\mathbf{q}$, we have
\begin{itemize}
\item 
$\es(h,k)\notin D$ for any $h$ and $k$ such that $h<\theta_{\mathbf{i},b}^{-1}(\pi(h))-2LN_4$ and $h\cdot\tsig(h)-k$ is even, and 
\item
$\es(h,k)\notin D$ for any $h$ and $k$ such that $h<\theta_{\mathbf{i},b}^{-1}(\pi(h))+2LN_4$ and $h\cdot\tsig(h)-k$ is odd.
\end{itemize} 
\end{lem}

\begin{NB}
Given $\sigma$, $\lambda$ and a monomial $\mathbf{q}$, 
there exists an integer  $B_2$ such that the following condition holds: 
for any $b>B_2$,

there exists an integer  $N_2$ such that the following condition holds: if $\theta\succ 3\cdot N_2$, 
then for any dimer configuration $D$ of type $(\sigma,\lambda,\nu,\theta)$ with weight $=\mathbf{q}$\,, we have
\begin{itemize}
\item 
$\es(h,k)\notin D$ for any $h$ and $k$ such that $h\in \pi^{-1}(j)$, $h<\theta^{-1}(j)-N_2\cdot L$ and $h\cdot\tsig(h)-k$ is even, and 
\item
$\es(h,k)\notin D$ for any $h$ and $k$ such that $h\in \pi^{-1}(j)$, $h<\theta^{-1}(j)+N_2\cdot L$ and $h\cdot\tsig(h)-k$ is odd.
\end{itemize} 
\end{lem}
\begin{proof}
First we take $\theta'\succ N_1$, where $N_1$ is an integer as in Lemma \ref{lem_52}.

We can take an integer $N'_2$ such that
\begin{itemize}
\item $e^{\mathrm{s}}(h,k)\notin D^{\mathrm{max}}_{\sigma,\lambda,\theta'}$ for any $h$ and $k$ such that $h\in \pi^{-1}(j)$, $h<\theta'^{-1}(j)-N'_2\cdot L$ and $h\cdot\tsig(h)-k$ is even, and 
\item $e^{\mathrm{s}}(h,k)\notin D^{\mathrm{max}}_{\sigma,\lambda,\theta'}$ for any $h$ and $k$ such that $h\in \pi^{-1}(j)$, $h>\theta'^{-1}(j)+N'_2\cdot L$ and $h\cdot\tsig(h)-k$ is odd.
\end{itemize} 
Since we have only a finite number of dimer configuration of type $(\sigma,\lambda,\nu,\theta')$ with weight $\mathbf{q}$ and each dimer configuration has only finite difference from $D^{\mathrm{max}}_{\sigma,\lambda,\nu,\theta'}$, we can take an integer $N''_2$ such that
\begin{itemize}
\item $e^{\mathrm{s}}(h,k)\notin D$ for any $h$ and $k$ such that $h\in \pi^{-1}(j)$, $h<\theta'^{-1}(j)-N''_2\cdot L$ and $h\cdot\tsig(h)-k$ is even, and 
\item $e^{\mathrm{s}}(h,k)\notin D$ for any $h$ and $k$ such that $h\in \pi^{-1}(j)$, $h>\theta'^{-1}(j)+N''_2\cdot L$ and $h\cdot\tsig(h)-k$ is odd.
\end{itemize} 
Take an integer $N_2$ such that $N>N''_2$ and $N_2\cdot \delta\succ\theta'$. 
For $\theta\succ N_2\cdot \delta$, we can take a minimal expression $\mathbf{i}$ such that $\theta_{\mathbb{i}}=\theta$ and $\theta_{\mathbf{i}_b}=\theta$ for some $b$. 
Then the claim follows from the fourth one in Lemma \ref{lem_52}.
\end{proof}
\end{NB}

\subsubsection{}
We assume that 
\[
\theta_{\mathbf{i},b}^{-1}(1/2)<\theta_{\mathbf{i},b}^{-1}(3/2)<
\cdots <\theta_{\mathbf{i},b}^{-1}(L-1/2)
\]
for any $b>0$.

Given $\sigma$, $\lambda$ and a monomial $\mathbf{q}$, 
take $B_5$ such that $B_5>b(2N_4)$ and $B_5>B_1$. 
The following theorem is the main result of this section:
\begin{thm}\label{thm_limit}
For any $b>B_5$, we have a bijection between the following two sets:
\begin{itemize}
\item dimer configurations of type $(\sigma,\lambda,\nu,\theta_{\mathbf{i},b})$ with weights $\mathbf{q}$,
\item data $(\vec{\Lambda},\vec{\nu})$ as in \S \ref{subsub_515} of type $(\sigma,\lambda,\nu)$ with weights $\mathbf{q}$.
\end{itemize}
\end{thm}
\begin{proof}
First, we divide the $(x,y)$-plane into the following $2L+1$ areas:
\begin{align*}
&C_{j}:=\{\theta^{-1}(j)-2LN_4<x<\theta^{-1}(j)+2LN_4\}\quad (j\in I_{\mathrm{h}}),\\
&C_{0}:=\{x<\theta^{-1}(1/2)-2LN_4\},\\
&C_{i}:=\{\theta^{-1}(i-1/2)+2LN_4<x<\theta^{-1}(i+1/2)-2LN_4\}\quad (1\leq i\leq L-1),\\
&C_{L}:=\{\theta^{-1}(L-1/2)+2LN_4<x\}.
\end{align*}
By Lemma \ref{lem_limit}, in the area $C_j$ we have
\begin{itemize}
\item $e^{\mathrm{s}}[h,k]\notin D$ for any $h$ and $k$ such that $\pi(h)>j$ and $h\cdot\tsig(h)-k$ is even, and 
\item $e^{\mathrm{s}}[h,k]\notin D$ for any $h$ and $k$ such that $\pi(h)<j$ and $h\cdot\tsig(h)-k$ is odd.
\end{itemize}
Removing these edges, we get a new graph.
A ``face" of the new graph is a union of $L$-tuple of elements in $\mathcal{F}$. Regarding such a union as a ``hexagon", the dimer configuration $D$ gives a dimer configuration for the hexagon lattice, in other words, a $3$-dimensional diagram. Let $\Lambda^{(j)}$ denote this $3$-dimensional diagram. (See Example \ref{ex4}.)

Similarly, in the area $C_j$ we have 
\begin{itemize}
\item $e^{\mathrm{s}}[h,k]\notin D$ for any $h$ and $k$ such that $\pi(h)>i$ and $h\cdot\tsig(h)-k$ is even, and 
\item $e^{\mathrm{s}}[h,k]\notin D$ for any $h$ and $k$ such that $\pi(h)<i$ and $h\cdot\tsig(h)-k$ is odd.
\end{itemize}
Removing these edges, we get a new graph, which is an infinite disjoint union of zigzag paths.
For each zigzag path, we have two choices of perfect matching and so the dimer configuration $D$ gives a Young diagram $\nu^{(i)}$.
We can verify that the datum $(\vec{\Lambda},\vec{\nu})$ satisfies the conditions in \S \ref{subsub_515}. Note that the reverse construction also works.

We have to check the correspondence above respects the weights.
Note that all edges of in the area $C_i$ have weights $=1$.
By \eqref{eq_weight3}, the contribution of the part in the area $C_j$ is given by
\[
\begin{cases}
w(\nu^{(j-1/2)};q^{(s_{\mathbf{i}}(j))}_+,Q_+)\cdot w(\nu^{(j+1/2)};{\bigl(q^{(s_{\mathbf{i}}(j))}_+\bigr)}^{-1}Q,Q_-)\cdot  w(\Lambda^{(j)}),& \sigma(j)=+,\\
w({}^{\mathrm{t}}\nu^{(j-1/2)};q^{(s_{\mathbf{i}}(j))}_+,Q_+)\cdot w({}^{\mathrm{t}}\nu^{(j+1/2)};{\bigl(q^{(s_{\mathbf{i}}(j))}_+\bigr)}^{-1}Q,Q_-)\cdot  w(\Lambda^{(j)}),& \sigma(j)=-.
\end{cases}
\]
Combining these contributions, we get the claim.

\begin{NB}
Given $\vec{\nu}$, we put 
\[
\Lambda^{(j)}_{\lambda,\nu,\vec{\nu}}:=\begin{cases}
\Lambda^{\mathrm{min}}_{\lambda^{(j)},\nu^{(j+1/2)},{}^{\mathrm{t}}\nu^{(j-1/2)}},& \sigma(j)=+,\\
\Lambda^{\mathrm{min}}_{\lambda^{(j)},{}^{\mathrm{t}}\nu^{(j-1/2)},\nu^{(j+1/2)}},& \sigma(j)=-,
\end{cases}
\]
and
\[
\vec{\Lambda}_{\lambda,\nu,\vec{\nu}}:=\bigl(\Lambda^{(j)}_{\lambda,\nu,\vec{\nu}}\bigr).
\]
Let $D_{\plmt}(\vec{\nu})$ be the dimer configuration of type $(\plmt)$ corresponding to the datum $\bigl(\vec{\Lambda}_{\lambda,\nu,\vec{\nu}},\vec{\nu}\bigr)$. 
Note that
\[
w\Bigl(\vec{\Lambda}_{\lambda,\nu,\vec{\emptyset}},\vec{\emptyset}\,\Bigr)
=
w_{\plt}\Bigl(D_{\plmt}\bigl(\vec{\emptyset}\,\bigr)\Bigr)
=1.
\]
First, we can verify
\[
w\bigl(\vec{\Lambda}_{\lambda,\nu,\vec{\nu}},\vec{\nu}\,\bigr)
=
w_{\plt}\bigl(D_{\plmt}(\vec{\nu}\,)\bigr)
\]
for any $\vec{\nu}$ by induction with respect to $\sum |\nu^{(i)}|$.
Then, we can check the claim for arbitrary datum by the induction with respect to 
\[
\sum_j |\Lambda^{(j)}\backslash \Lambda^{(j)}_{\lambda,\nu,\vec{\nu}}|.
\]
\end{NB}
\end{proof}
\begin{ex}\label{ex4}
We take $\sigma$ as in Example \ref{ex1} and $\lambda=\emptyset$.
Assume that $\theta(1/2)=N+1/2$ and $\theta(5/2)=-N+5/2$ for $N\gg 0$.
In Figure \ref{fig23}, we show the weight \textup{(}after putting $q_+=q_-=q_0^{1/2}$\,\textup{)} of edges in the area $C_{1/2}$. 
We can idenfity the graph in the area $C_{1/2}$ with a hexagon lattice as in Figure \ref{fig23}.
\begin{figure}[htbp]
  \centering
  \input{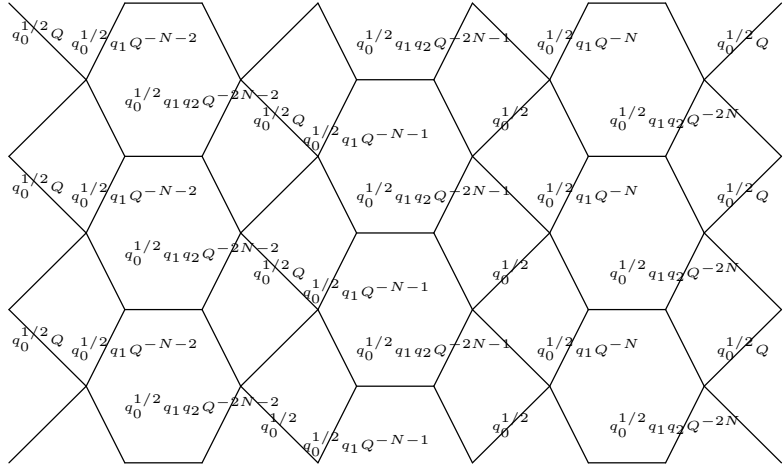}
  \caption{the graph in the area $C_{1/2}$}
\label{fig23}
\end{figure}
\begin{figure}[htbp]
  \centering
  \input{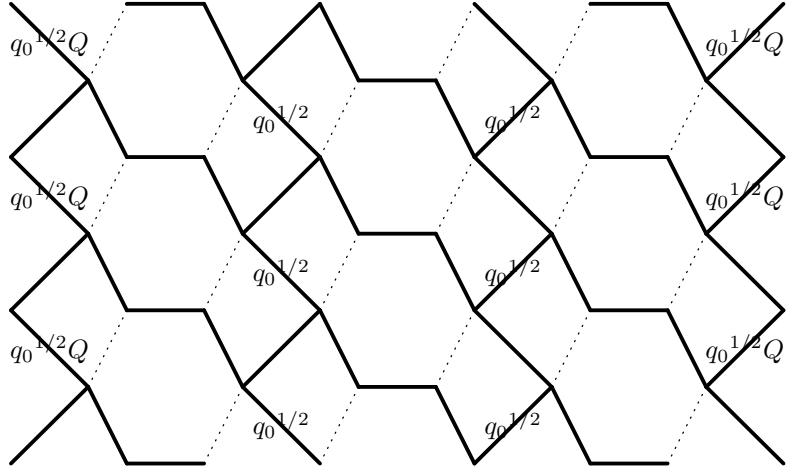}
  \caption{identification with a hexagon lattice}
\label{fig24}
\end{figure}
\end{ex}
\begin{rem}
In general, we have the permutation  $s_{\mathbf{i}}\in \mathfrak{S}_{\Ih}$ of the set $\Ih$ satisfying the following condition: for sufficiently large $b$ we have
\[
\theta_{\mathbf{i},b}^{-1}(s_{\mathbf{i}}(1/2))<\theta_{\mathbf{i},b}^{-1}(s_{\mathbf{i}}(3/2))<
\cdots <\theta_{\mathbf{i},b}^{-1}(s_{\mathbf{i}}(L-1/2)).
\]
The permutation $s_{\mathbf{i}}$ determines the direnction in which we take limit in the space of stability conditions.
It is the refine topological vertex associated to $Y_{\sigma\circ s_{\mathbf{i}}}$ what we get in the limit.
\end{rem}

\subsection{Conclusion}
Note that we have
\[
\bigcup_{b=1}^\infty\Lambda_{\theta_{\mathbf{i},b}}^{\mathrm{re},+}=\Lambda_+^{\mathrm{re},+}.
\]
Combining Theorem \ref{thm_410} and Theorem \ref{thm_limit}, we have the following:
\begin{thm}
\begin{align*}
\mathcal{Z}^{\mathrm{RTV}}_{\sigma,\lambda,\nu}
&=
\mathcal{Z}_{\sigma,\lambda,\nu}^{\mathrm{NCDT}}
\cdot
\prod_{\alpha\in\Lambda_+^{\mathrm{re},+}}\Biggr(\prod_{(h,h')\in B^{\alpha,+}_{\sigma,\lambda}}(1-\sigma(\alpha)\cdot w_\lambda(h')/w_\lambda(h))^{\sigma(\alpha)}\\
&\hspace{50mm}\prod_{(h,h')\in B^{\alpha,-}_{\sigma,\lambda}}(1-\sigma(\alpha)\cdot w_\lambda(h')/w_\lambda(h))^{-\sigma(\alpha)}\Biggl).
\end{align*}
\textup{(}See \eqref{eq_re++}, \eqref{eq_Bpm1} and \eqref{eq_sigma_alpha} for the notations.\textup{)}
\end{thm}
Since the second term in the right-hand side does not depend on $\nu$, we have the following:
\begin{cor}\label{cor_56}
\[
\mathcal{Z}^{\mathrm{RTV}}_{\sigma,\lambda,\nu}\big/\mathcal{Z}_{\sigma,\lambda,\vec{\emptyset}}^{\mathrm{RTV}}=
\mathcal{Z}_{\sigma,\lambda,\nu}^{\mathrm{NCDT}}\big/\mathcal{Z}_{\sigma,\lambda,\vec{\emptyset}}^{\mathrm{NCDT}}.
\]
\end{cor}
Combining Theorem \ref{thm_412} and Theorem \ref{thm_limit}, we have the following:
\begin{thm}
\begin{align*}
&\mathcal{Z}_{\sigma,\lambda,\nu}^{\mathrm{RTV}}\Big|_{q_+=q_-=(q_0)^{1/2}}=\\
&\quad \mathcal{Z}_{\sigma,\lambda,\nu}^{\mathrm{NCDT}}\Big|_{q_+=q_-=(q_0)^{1/2}}
\cdot
\prod_{
\alpha\in \Lambda_+^{\mathrm{re},+}}
\bigl(1-\sigma(\alpha)\cdot q^{\alpha}\bigr)^{\sigma(\alpha)\cdot \bigl(\alpha^0+c_\lambda[j_-(\alpha)]-c_\lambda[j_+(\alpha)]\bigr)}.
\end{align*}
\textup{(}See \eqref{eq_re++}, \eqref{eq_Bpm1} and \eqref{eq_sigma_alpha} for the notations.\textup{)}
\end{thm}
Since the second term in the right-hand side depend only on $c_\lambda[j]$'s but not on $\lambda$ and $\nu$, we have the following:
\begin{cor}\label{cor_58}
If $c_\lambda[j]=0$ for any $j$, then we have
\[
\mathcal{Z}_{\sigma,\lambda,\nu}^{\mathrm{RTV}}\big/\mathcal{Z}^{\mathrm{RTV}}_{\sigma,\vec{\emptyset},\vec{\emptyset}}\,\Big|_{q_+=q_-}
=
\mathcal{Z}_{\sigma,\lambda,\nu}^{\mathrm{NCDT}}\big/\mathcal{Z}_{\sigma,\vec{\emptyset},\vec{\emptyset}}^{\mathrm{NCDT}}\,\Big|_{q_+=q_-}
.
\]
\end{cor}

\bibliographystyle{amsalpha}
\bibliography{bib-ver6}

\end{document}